\documentclass{article}
\usepackage[utf8]{inputenc}

\makeatletter
\renewcommand*\l@section{\@dottedtocline{1}{0em}{2em}}
\renewcommand*\l@subsection{\@dottedtocline{2}{1.5em}{3em}}
\makeatother

\usepackage{amsmath,amssymb,amsfonts}
\usepackage{physics}
\usepackage{dsfont}
\usepackage{graphicx}

\usepackage{color}
\definecolor{purple}{rgb}{0.5,0,0.5}

\usepackage{subcaption}

\definecolor{cof}{RGB}{219,144,71}
\definecolor{pur}{RGB}{186,146,162}
\definecolor{greeo}{RGB}{91,173,69}
\definecolor{greet}{RGB}{52,111,72}

\newcommand{\ra}{\rightarrow}

\newcommand{\End}{{\rm End}}

\newcommand{\NN}{{\mathbb N}}
\newcommand{\CC}{{\mathbb C}}
\newcommand{\ZZ}{{\mathbb Z}}
\newcommand{\RR}{{\mathbb R}}

\newcommand{\cA}{\mathcal A}
\newcommand{\cB}{\mathcal B}

\newcommand{\cZ}{\mathcal Z}

\newcommand{\cM}{\mathcal M}

\newcommand{\pinm}{\text{pin}^-}

\newcommand{\fotimes}{\mathbin{\widehat\otimes}}
\newcommand{\bigfotimes}{\mathbin{\widehat\bigotimes}}

\newcounter{save}
\newenvironment{AxAlign}
   {\setcounter{save}{\theequation}
     \setcounter{equation}{0}%
     \renewcommand\theequation{A\arabic{equation}}
     \align}
   {\endalign\setcounter{equation}{\thesave}}

\usepackage[left=0.95in,right=0.95in,top=0.85in,bottom=0.85in]{geometry}

\usepackage{hyperref}

\hypersetup{colorlinks=true,urlcolor=[rgb]{0,0,0.5},citecolor=[rgb]{0,0.5,0},linkcolor=[rgb]{0,0.5,0}}

\renewcommand{\thesection}{\Roman{section}}

\title{Diagrammatic State Sums for 2D Pin-Minus TQFTs}
\author{Alex Turzillo\\
{\it California Institute of Technology, Pasadena, CA 91125}}

\begin{document}
\maketitle

\begin{abstract}
    The two dimensional state sum models of Barrett and Tavares are extended to unoriented spacetimes. The input to the construction is an algebraic structure dubbed half twist algebras, a class of examples of which is real separable superalgebras with a continuous parameter. The construction generates pin-minus TQFTs, including the root invertible theory with partition function the Arf-Brown-Kervaire invariant. Decomposability, the stacking law, and Morita invariance of the construction are discussed.
\end{abstract}

\tableofcontents

\pagebreak

\section{Introduction}

State sum constructions of quantum field theories extend Feynman's formulation of the time-sliced quantum mechanical path integral to theories of positive spatial dimension. They are closely related to lattice models, which are expected to generate all consistent\footnote{Free of anomalies, such as the framing anomalies suffered by Reshetikhin-Turaev theories with nonzero chiral central charge.} quantum field theories by a continuum limit. In the case of topological theories, which are sensitive only to the spacetime topology (rather than a metric), the study of state sums has been particularly fruitful, with applications in mathematics --- perhaps most famously, to knot theory \cite{K87} --- as well as in physics. Merits of the state sum approach include that its algebraic input is simpler than the continuum data and that quantities of interest may be computed by local algorithms. This, however, comes at the cost of redundancy, as lattice realizations are not unique. As we will see, this trade-off essentially reflects the difference between certain algebraic structures and their Morita classes.

Topological quantum field theories (TQFTs) have recently gained prominence in condensed matter physics due to their connection to topological phases of matter. It is claimed that the field theories encode the universal, long-distance effective behavior --- the ``phase'' --- of gapped quantum systems, which means characterizing their responses to topological probes and reproducing the ground state expectation values of nonlocal order parameters \cite{SR16,SSGR17}. State sum constructions of the field theories are directly related to the gapped lattice models that live at renormalization group fixed points \cite{KTY16a,SR16,KTY16b}. According to this picture, a sensitivity of the theory to a spin structure, in addition to topology, captures the response of massive fermions in the gapped system to boundary conditions. Such field theories are known as spin-TQFTs. When a gapped system has a time-reversal symmetry, its effective field theory is insensitive to the orientation of spacetime and is defined on all unoriented spacetimes.\footnote{The path integral on nonorientable spacetimes computes the time-reversal symmetry protected trivial (SPT) order \protect\cite{SR16,SSGR17}.} When fermions transform under time-reversal symmetry with $T^2=\mp 1$, the appropriate geometric structure is a $\text{pin}^\pm$ structure. Of particular physical relevance are $\text{pin}^-$ theories in two (spacetime) dimensions and their relationship with time-reversal-invariant Majorana chains, which have been known for some time to have an interesting interacting gapped phase classification \cite{FK11}.

Given the usefulness of state sum models for purely topological theories, it is natural to ask whether spin- and $\text{pin}^\pm$-TQFTs yield state sums as well. The case of spin theories in two spacetime dimensions was recently studied by Barrett and Tavares \cite{BT14} (see ref. \cite{NR15} for an alternate approach). They exploit the relation between spin structures on a surface $M$ and immersions of $M$ into $\RR^3$ to construct, for each spin surface, a ribbon diagram, the twists and crossings of which keep track of the spin structure. Their state sum models are then computed locally on this discrete realization of the spin geometry.

The main result of our paper is a state sum construction for two dimensional $\text{pin}^-$ theories. Our approach extends that of ref. \cite{BT14} to unoriented spacetimes. The state sums amount to discretizations of all unitary invertible (as well as many non-unitary and/or non-invertible) field theories with this structure, in particular the Arf-Brown-Kervaire theory, which was recently studied along with its connection to Majorana chains in ref. \cite{DG18}. A broad class of them has a simple algebraic characterization in terms of certain real superalgebras. From this perspective, the eight distinct powers of the Arf-Brown-Kervaire theory (the eight phases of time-reversal-invariant Majorana chains) arise from the eight Morita classes of central simple real superalgebras,\footnote{As discussed below, the state sums for the non-central algebras describe the two symmetry-broken theories.} a connection which has been noted previously in the context of tensor network states \cite{FK11,TY17,BWHV17}. In topological theories, the state sum data has an interpretation as the space of states on the interval \cite{LP06}; similarly, the real Clifford algebras $C\ell_{n,0}\RR$, $n=0,\ldots,7$, whose state sums are the eight invertible $\text{pin}^-$ theories, have to do with Majorana zero modes localized at the endpoints of the open chain.

The structure of the paper is as follows. In section \ref{pin2D}, we review some elementary facts about pin structures on closed surfaces and cobordisms and their relation to codimension one immersions and quadratic enhancements. Diffeomorphism classes of these objects and their classification by the Arf-Brown-Kervaire invariant are discussed. We also derive a simple expression for the evaluation of the quadratic enhancement on an embedded curve in terms of its ribbon diagram. In section \ref{ribht}, we show how to construct a ribbon diagram from an immersed surface and evaluate its state sum. Imposing invariance under re-triangulation and regular homotopy, we derive the defining axioms of a half twist algebra. The state spaces of the associated $\text{pin}^-$-TQFT are constructed as well. In section \ref{supalgabk}, we specialize to a class of half twist algebras related to real superalgebras. Decomposability and stacking are understood on the level of these algebras, and it is shown that Morita equivalent algebras define the same theory. We explicitly compute the path integrals for the Euler and Arf-Brown-Kervaire theories and discuss the classification of invertible $\text{pin}^-$-TQFTs.

\section{Pin geometry in two dimensions}\label{pin2D}

\subsection{Pin structures, immersions, and quadratic enhancements}

The goal of this section is to review the following equivalences:
$$\left\{\begin{array}{c}\text{pin}^-\text{ structures / isom.}\\/\\\text{pin}^-\text{-diffeomorphism}\\=\\\text{pin}^-\text{-diffeo. classes}\end{array}\right\}\begin{array}{c}\leftrightarrow\\\\\leftrightarrow\\\\\leftrightarrow\end{array}\left\{\begin{array}{c}\text{quadratic enhancements}\\/\\\text{lin. aut. with }q'=q\circ\alpha\\=\\\text{quadratic enh. / equiv.}\end{array}\right\}\begin{array}{c}\leftrightarrow\\\\\leftrightarrow\\\\\leftrightarrow\end{array}\left\{\begin{array}{c}\text{immersions / reg. homot.}\\/\\\text{diffeo. with }f=g\circ\phi\\=\\\text{imm. surf. / reg. homot.}\end{array}\right\}$$

Pin structures generalize spin structures to unoriented smooth manifolds. The structure group $O(n)$\footnote{A Riemannian metric is required to reduce the structure group from $GL_n\RR$ to $O(n)$.} of an unoriented manifold has two double covers $Pin^-(n)$ and $Pin^+(n)$, which differ in the behavior of the lifts $\tilde r$ of odd reflections $r\in O(n)$: in $Pin^\pm(n)$, they square to $\tilde r^2=\pm 1$ \cite{ABS64}. A $\text{pin}^\pm$ structure on an unoriented manifold is a principal $Pin^\pm(n)$ bundle with a $2$-fold covering of the orthogonal frame bundle that restricts to the double cover $\rho:Pin^\pm(n)\rightarrow O(n)$ on fibers. The following discussion of $\text{pin}^\pm$ structures is adopted from ref. \cite{KT89}. In terms of an open cover on $M$, it is a global lift of the $O(n)$-valued transition functions $t_{ij}$ to $s_{ij}\in Pin^\pm(n)$. The triple overlap condition $t_{ij}t_{jk}t_{ki}=1$ ensures that any local lifts $\rho:s_{ij}\mapsto t_{ij}$ satisfy $s_{ij}s_{jk}s_{ki}=o_{ijk}\in\ker\rho\simeq\ZZ/2$. By looking at the quadruple overlap, one sees that the signs $o_{ijk}$ form a \v Cech $2$-cocycle. Local lifts are acted on transitively by $\ker\rho$-valued $1$-cochains $A$ as $s_{ij}\mapsto s_{ij}A_{ij}$, which shifts $o$ by the coboundary $\delta A$. The class $[o]\in H^2(M;\ZZ/2)$ is the obstruction to a global lift, or $\text{pin}^\pm$ structure, and is $w_2+w_1^2$ for $\pinm$ and $w_2$ for $\text{pin}^+$, where the $w_i$ denote the Stiefel-Whitney classes of the tangent bundle of $M$. Two $\text{pin}^\pm$ structures are regarded as isomorphic if they are related by a transformation $s_{ij}\mapsto \lambda_is_{ij}(\lambda_j)^{-1}$, $\lambda_i\in Pin^\pm(n)$. If $A$ is closed\footnote{\v Cech cocycles $A\in Z^1(M;\ZZ/2)$ are often referred to as $\ZZ/2$-gauge fields.} and $s$ is a $\text{pin}^\pm$ structure, the lift $sA$ is again a $\text{pin}^\pm$ structure, and the two are isomorphic iff $A$ is a coboundary $\delta\lambda$; thus, assuming $[o]$ vanishes, isomorphism classes of $\text{pin}^\pm$ structures on $M$ form a torsor for $H^1(M;\ZZ/2)$. Our focus will be on surfaces and their $\pinm$ structures, or simply ``pin structures.'' The obstruction class vanishes in two dimensions, so each surface supports exactly $|H^1(M;\ZZ/2)|$ pin structures, up to isomorphism.

Another characterization of pin structures on a surface $M$ can be given in terms of immersions of $M$ into $\RR^3$. Two immersions are said to be \emph{regular homotopic} if they are connected by a smooth $1$-parameter family of immersions \cite{P84}. Immersions of a surface $M$ into $\RR^3$ fall into $|H^1(M;\ZZ/2)|$ regular homotopy classes \cite{P84,JT66,S77}, one for each isomorphism class of pin structure on $M$. The pin structure corresponding to an immersion is obtained by pulling back the standard pin structure on $\RR^3$ by the immersion \cite{KT89}. Two immersions $f, g$ are \emph{equivalent} if there exists a diffeomorphism $\phi$ of $M$ such that $f=g\circ\phi$, and these equivalence classes, called \emph{immersed surfaces}, are said to be regular homotopic if their representative immersions are \cite{P84}. Equivalence of immersions corresponds to pin diffeomorphism of the corresponding pin surfaces.

Pin structures on surfaces have a third characterization: their isomorphism classes are in bijective correspondence with \emph{quadratic enhancements of the intersection form} \cite{KT89}; that is, functions
\begin{equation}
    q:H_1(M;\ZZ/2)\rightarrow\ZZ/4
\end{equation}
such that
\begin{equation}\label{quad}
    q(x+y)=q(x)+q(y)+2\cdot\langle x,y\rangle,
\end{equation}
where $2\cdot$ embeds $\ZZ/2$ into $\ZZ/4$ as a subgroup and $\langle\cdot,\cdot\rangle$ denotes the intersection form on $M$. In ref. \cite{KT89} Kirby and Taylor demonstrate how to build a quadratic enhancement from a pin structure, while in ref. \cite{P84} Pinkall does the same from its associated immersion. Since the constructions are similar, below we will focus solely on the latter. Every quadratic enhancement arises from both a pin structure and an immersion, and the constructions are isomorphism and regular homotopy invariant, respectively. We say that two quadratic enhancements $q,q'$ are \emph{equivalent} if they are related as $q'=q\circ\alpha$ by a linear automorphism $\alpha$ of $H_1(M;\ZZ/2)$. As all linear automorphisms $\alpha$ that preserve the intersection form are induced by diffeomorphisms of $M$ \cite{P84,MP78}, all equivalences of quadratic enhancements arise from equivalences of immersions. A pin diffeomorphism that covers a diffeomorphism $\phi$ of the base space $M$ induces an equivalence $q'=q\circ\phi_*$ on the associated quadratic forms. Quadratic enhancements form a torsor for $H^1(M;\ZZ/2)$ by the action $q\mapsto q+2\cdot A$, with respect to which the correspondence with pin structures is equivariant \cite{KT89}.

\subsection{The quadratic enhancement as a self-linking number}\label{quadsubsec}

Let us now follow ref. \cite{P84} in constructing a quadratic enhancement from an immersion. Begin by defining a function $\tilde q_f$ that takes closed loops in $M$ to their self-linking numbers. To be precise, $\tilde q_f$ is defined on smooth embeddings $\gamma:S^1\rightarrow M$ such that $f\circ\gamma:S^1\rightarrow\RR^3$ is also an embedding. Images of such embeddings have embedded tubular neighborhoods (``ribbons'') $N_\gamma$. The self-linking number is given by the linking number of the loop $f\circ\gamma$ with the loop obtained by pushing $f\circ\gamma$ along $N_\gamma$:
\begin{equation}
    \tilde q_f(\gamma)=\text{link}(f\circ\gamma,f(\partial N_\gamma)).
\end{equation}
Under regular homotopy, $\tilde q_f$ is stable only modulo $4$; moreover, it depends only on the $\ZZ/2$-homology class $[\gamma]\in H_1(M;\ZZ/2)$ and defines a map $q_f$ on $H_1(M;\ZZ/2)$ satisfying the quadratic enhancement condition \eqref{quad}.

\begin{figure}
    
    \centering
	\includegraphics[width=15cm]{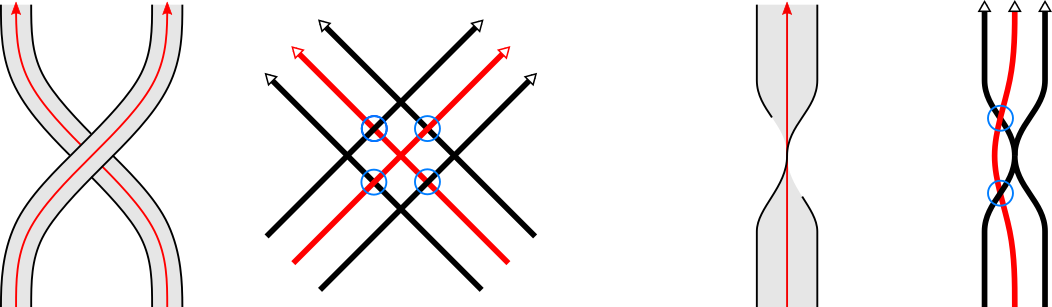}

  \caption{Following Pinkall \protect\cite{P84}, give the core (red) and edges (black) of the ribbon a particular orientation. Then compute the linking number of the red lines with the black lines. The crossing has four red-black intersections, all of the same parity. The half twist has two red-black intersections of the same parity.}
  \label{selflinking}
\end{figure}

By projecting a ribbon onto $\RR^2$ and obtaining a ribbon diagram, its self-linking number may be computed by a local algorithm. As is discussed in greater detail in section \ref{ribbon}, one may use regular homotopy so that the projection $\RR^3\ra\RR^2$ onto the $xy$-plane is an immersion of $N_\gamma$ at all but finitely many points where the ribbon makes a half twist (left or right handed). The image of the curve $\gamma$ may be taken to cross itself transversely and away from these points. Away from the twists and crossings, the self-linking number is zero. As demonstrated in Figure \ref{selflinking}, each right handed half twist contributes $+1$ to $\tilde q_f(\gamma)$; likewise, each left handed half twist contributes $-1$. Each crossing contributes $\pm 2$. In total,
\begin{equation}\label{quadlocal}
    \tilde q_f(\gamma)=(\#\text{ r.h. twists})-(\#\text{ l.h. twists})+2\cdot(\#\text{ crossings})\mod 4.
\end{equation}

\subsection{The Arf-Brown-Kervaire invariant}

The Arf-Brown-Kervaire (ABK) invariant of a surface $M$ with quadratic enhancement $q$ is defined as
\begin{equation}\label{abk}
    \text{ABK}(M,q)=\frac{1}{\sqrt{|H_1(M;\ZZ/2)|}}\sum_{x\in H_1(M;\ZZ/2)}e^{i\pi q(x)/2}.
\end{equation}
It is valued in eighth roots of unity and has the nice property that two quadratic enhancements on $M$ have the same ABK invariant if and only if they are equivalent \cite{B72}. In other words, the ABK invariant is well-defined on diffeomorphism classes of pin surfaces as well as on immersed surfaces. The ABK invariant determines the pin bordism class of the pin surface and so defines an isomorphism $\Omega_2^\text{pin}(\text{pt})\xrightarrow{\sim}\ZZ/8$.

\subsection{Decomposition of pin surfaces}\label{direct}

Every closed unoriented surface may be decomposed as a connect sum of tori and real projective planes. Each of these building blocks has two diffeomorphism classes of pin structures. On the torus, there are four isomorphism classes of pin structures given by a choice of NS (bounding, antiperiodic) or R (non-bounding, periodic) boundary conditions around each independent $1$-cycle. Pin diffeomorphisms covering Dehn twists relate the NS-NS, NS-R, and R-NS classes. To see this, note that a Dehn twist induces a map $\{x',y'\}=\{x,x+y\}$ on a basis of $H_1(T^2;\ZZ/2)=\ZZ/2\times\ZZ/2$. Then use the rule \eqref{quad}: the NS-NS pin structure $q(x)=0, q(y)=0$ becomes the NS-R pin structure
\begin{equation}
    q(x')=q(x)=0,\qquad q(y')=q(x+y)=q(x)+q(y)+2\cdot\langle x,y\rangle=2.
\end{equation}
These pin structures are not diffeomorphic to the R-R pin structure. One may also use \eqref{abk} to see that the NS-NS, NS-R, and R-NS pin structures have ABK invariant $+1$ (and so are diffeomorphic to each other), while the R-R pin structure has ABK invariant $-1$. Moreover, since the ABK invariant determines the bordism class, this calculation shows that the NS-NS pin structure bounds a solid torus, while the R-R pin structure is non-bounding. On the real projective plane, there are two isomorphism classes of pin structure. To see this, note that $H_1(\RR P^2;\ZZ/2)=\ZZ/2$, the generator $z$ of which is represented by $1$-sided (i.e. orientation-reversing) curve and has self-intersection $\langle z,z\rangle=1$. Since $q(0)=0$, the rule \eqref{quad} says
\begin{equation}
    0=q(z)+q(z)+2\cdot\langle z,z\rangle=2q(z)+2\mod 4,
\end{equation}
so there are two isomorphism classes of pin structures given by $q(z)=1$ and $q(z)=3$. These are non-diffeomorphic since they have ABK invariants $\exp(i\pi/4)$ and $\exp(7i\pi/4)$, respectively. Call them $\RR P^2_1$ and $\RR P^2_7$.

The pin structures on other surfaces may be readily understood from their connect sum decompositions. For example, the Klein bottle decomposes as $K\simeq\RR P^2\#\RR P^2$. Let $z_1,z_2$ denote the generating $1$-(co)cycles of the real projective planes. In this basis, the four quadratic enhancements are $q=(1,1), (1,3), (3,1), (3,3)$. In the familiar basis of $H_1(K;\ZZ/2)$ given by the orientation-preserving curve $x=z_1+z_2$ and orientation-reversing curve $y=z_2$, the possibilities are $q=(2,1), (0,3), (0,1), (2,3)$. They have ABK invariants $+i$, $+1$, $+1$, and $-i$, so there are three diffeomorphism classes of pin structures on $K$, one of which is null-bordant.

\subsection{Pin bordism and TQFT}\label{pincob}

Our discussion so far has focused on closed surfaces. To define pin TQFTs, it is necessary to also understand pin one manifolds and the bordisms between them. There are two connected one dimensional pin manifolds given by the NS and R spin structures on the circle. A pin manifold with boundary induces a pin structure on its boundary, and a pin bordism between pin one manifolds $S_0$ and $S_1$ is a pin surface $M$ whose boundary, with induced pin structure, is $S_0\sqcup S_1$.

\begin{figure}

	\centering
	\includegraphics[width=12.5cm]{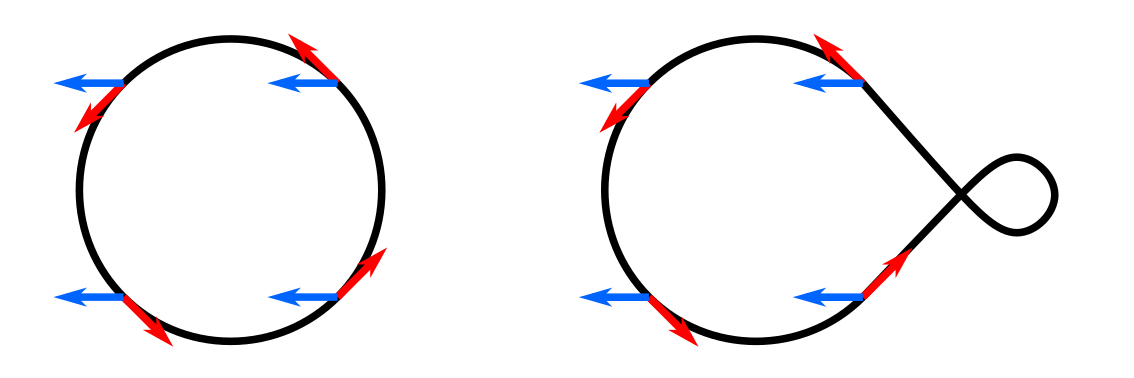}

  \captionof{figure}{Two examples of immersions of the circle in the plane, with turning numbers $1$ (left) and $0$ (right) defined as the winding of a tangent frame (red) relative to a constant vector field (blue). This number mod $2$ determines the induced (s)pin structure on the circle: NS for odd, R for even.}
  \label{circles}
\end{figure}

Each of the two pin structures on the circle is related to a class of immersed circles in the plane, depicted in Figure \ref{circles}. Fix two planes $\RR^2_0, \RR^2_1$ normal to the $y$-axis. An immersion of the cobordism $(S_0,S_1,M)$ is an immersion of $M$ such that $S_0, S_1$ lie in $\RR^2_0, \RR^2_1$, respectively. A regular homotopy of the immersions of the cobordism is again a $1$-parameter family of immersions. We emphasize that at each value of the parameter, the boundaries $S_0, S_1$ are pinned to the planes $\RR^2_0, \RR^2_1$.

The theory of quadratic enhancements associated to pin surfaces with boundary requires more care than we will give it here. The idea is to extend the discussion of ref. \cite{S04}. Choose a set of basepoints $\partial_0M$ --- one on each connected component of $\partial M$, and let a pin structure on $(M,\partial_0M)$ be a pin structure on $M$ together with a trivialization of the $Pin^-(1)=\ZZ/4$ bundle over $\partial_0M$. Such pin structures should be (non-canonically) identified with quadratic enhancements of the intersection form on $H_1(M_*;\ZZ/2)\simeq H_1(M,\partial_0M;\ZZ/2)$, where $M_*$ is a closed pin surface obtained by sewing a punctured sphere into $M$.

A pin TQFT assigns state spaces $\cA_{NS}, \cA_R$ to the circles $S^1_{NS}, S^1_R$ and linear maps to the pin bordisms between them. In particular, the mapping cylinders associated to elements of the pin mapping class group of the circles defines a supervector space structure on the state spaces. A complete algebraic characterization of pin TQFTs would resemble the discussions of refs. \cite{MS06,TT06,DG18}. We will not give one here; instead our focus will be on the pin TQFTs that arise from the diagrammatic state sum construction introduced below.

\section{Ribbon diagrams and half twist algebras}\label{ribht}

A state sum model provides a combinatorial description of a theory like a TQFT or, in the present case, a pin TQFT. Focusing first on defining partition functions of closed spacetime manifolds, the idea is to define an invariant of discretized spacetimes, given as a weighted sum over colorings of a discretization. The weight assigned to a coloring is computed ``locally'' from contributions of local elements of the discretization. The requirement that the invariant is independent of the discretization imposes structure on the weights.

For example, in ref. \cite{FHK94} Fukuma, Hosono, and Kawai study two-dimensional \emph{topological} state sums, which are defined on triangulated surfaces and whose weights receive contributions from the faces and edges of the triangulation. Topological invariance --- that is, lack of dependence on the triangulation --- imposes Pachner move conditions on this algebraic data. The result is that the local tensors assigned to faces and edges form a separable algebra.

State sum models for pin TQFTs have a similar logic. A discretization of a pin surface is a triangulation together with an additional combinatorial structure representing a pin structure. Finding these structures and the equivalence relations under which they represent the same continuum structure is not easy. One approach is to find a local combinatorial structure, or \emph{marking}, as Novak and Runkel do for spin structures in ref. \cite{NR15}. This paper follows a different path, one based on the connection between pin structures and immersions into $\RR^3$. In the following, a discretization is a triangulation together with a choice of immersed surface. The construction is automatically invariant under equivalence of immersion, whereas invariance under regular homotopy is enforced by hand. The weights are products of tensors assigned to elements of the discretization. The requirement of invariance under change of discretization (Pachner moves and regular homotopy) means that these tensors satisfy several relations. The resulting algebraic structure is what we dub a \emph{half twist algebra} and extends the separable algebras of ref. \cite{FHK94} to allow for the theory's sensitivity to pin structure.

\subsection{Ribbon diagrams}\label{ribbon}

\begin{figure}

	\centering
	\includegraphics[width=15.7cm]{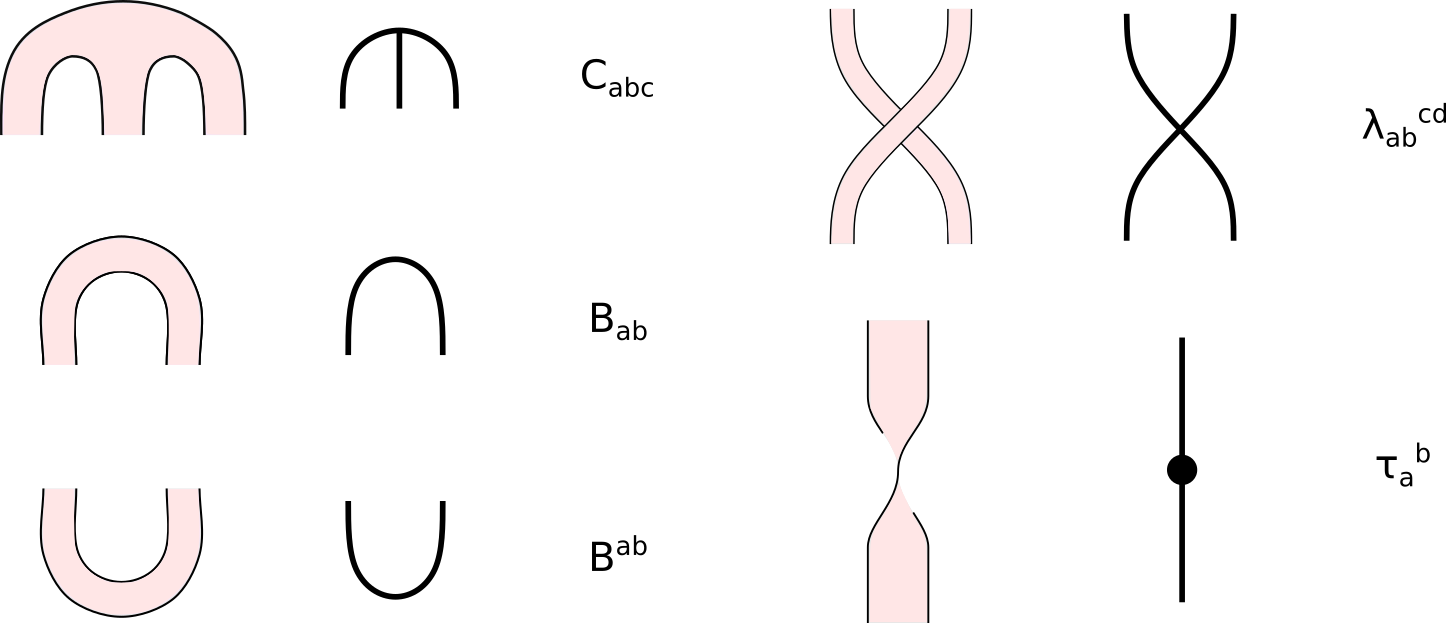}

  \caption{The five building blocks of ribbon diagrams satisfying the regularity conditions.}
  \label{blocks}
\end{figure}

We now construct a ribbon diagram from a triangulation of an immersed surface. Dual to the triangulation of the surface is a graph, which may be enlarged to a \emph{ribbon graph} by taking a regular neighborhood, the compliment of which in $M$ is one or more disks. Any immersion of $M$ is regular homotopic to one that is an embedding on the ribbon graph \cite{P84}. This embedded ribbon graph is passed through the projection $p:\RR^3\rightarrow\RR^2$ onto the $xy$-plane.\footnote{The ribbon diagrams associated to any two projections are related by rotation of the immersed surface in $\RR^3$, which is a regular homotopy. Since the state sum is, by construction, regular homotopy invariant, the choice of $p$ does not matter.} By regular homotopy, the projection can be made to satisfy certain regularity conditions. First, the projection is an immersion of the ribbon graph at all but finitely many points where the ribbon makes a \emph{half twist} \cite{K90}. Second, the edges of the graph intersect transversely in the image of $p$. Third, the graph is parallel to the $x$-direction at only finitely many ``critical points'' (nodes, caps, cups) where either all legs exit above the $x$-parallel or all legs exit below (no saddle points). Fourth, each node of the graph is located at a critical point with its three legs exiting below. Fifth, at most one of the following can occur at any point: a half twist, a crossing, and a critical point. In addition to the image of the projection, the helicities of the half twists (right or left handed) are recorded. Unlike diagrams typical in knot theory, ours do not record whether one strand crosses over or under the other at a crossing, as these two configurations are related by regular homotopy. A ribbon diagram satisfying the regularity conditions is composed of the five building blocks --- nodes, caps, cups, crossings, and half twists --- depicted in Figure \ref{blocks}.

\begin{figure}
\centering
\begin{subfigure}{.38\textwidth}
  \centering
  \includegraphics[width=\linewidth]{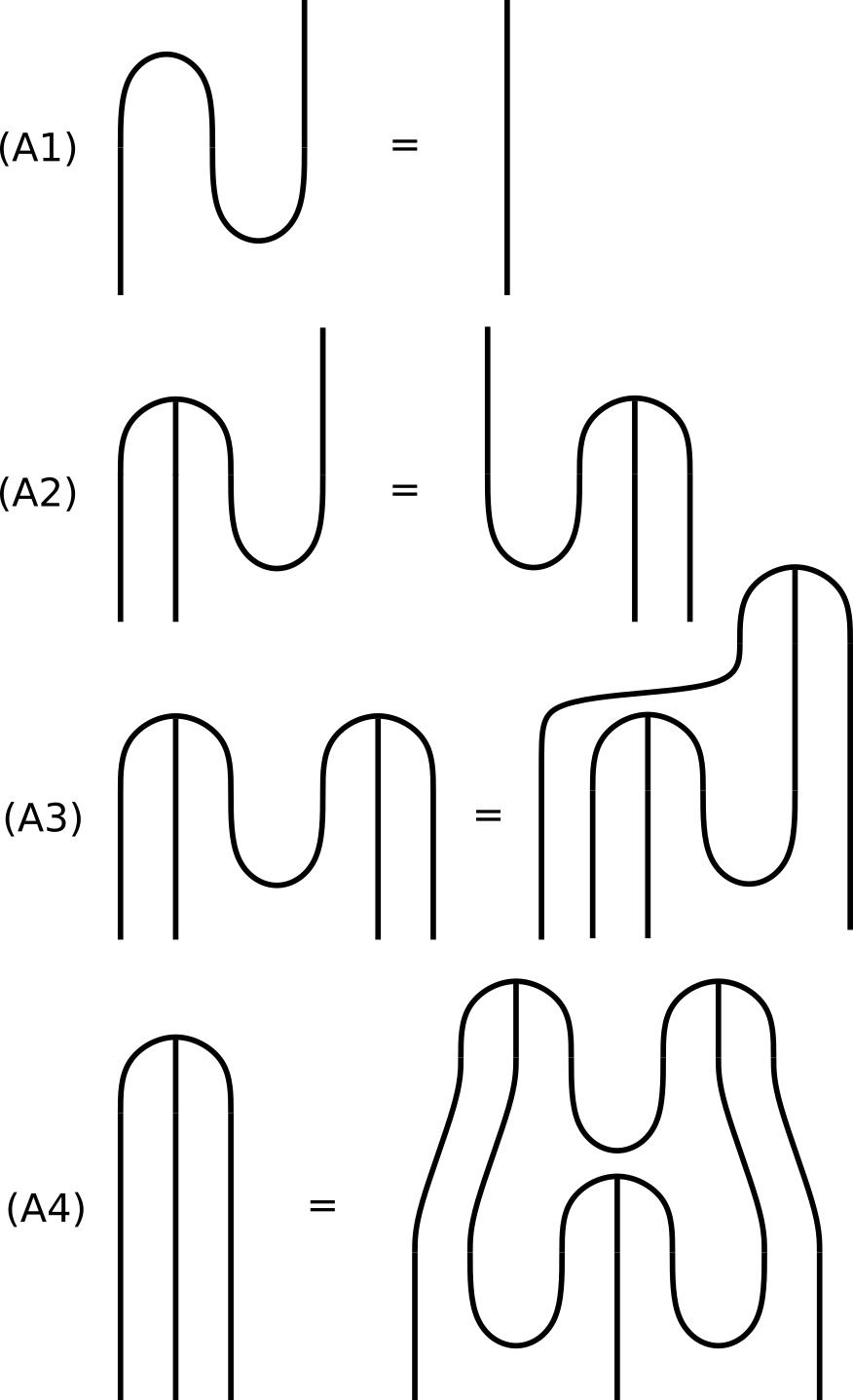}
\end{subfigure}%
\begin{subfigure}{.3\textwidth}
  \centering
  \includegraphics[width=0.8\linewidth]{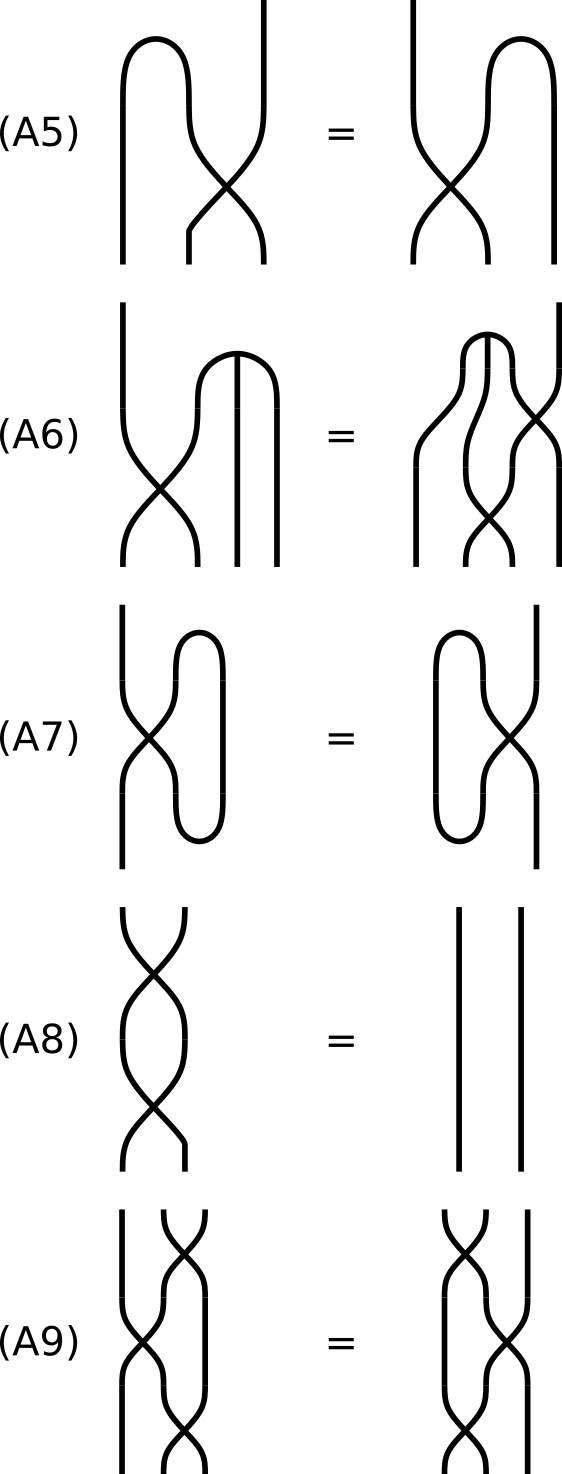}
\end{subfigure}
\begin{subfigure}{.3\textwidth}
  \centering
  \includegraphics[width=0.84\linewidth]{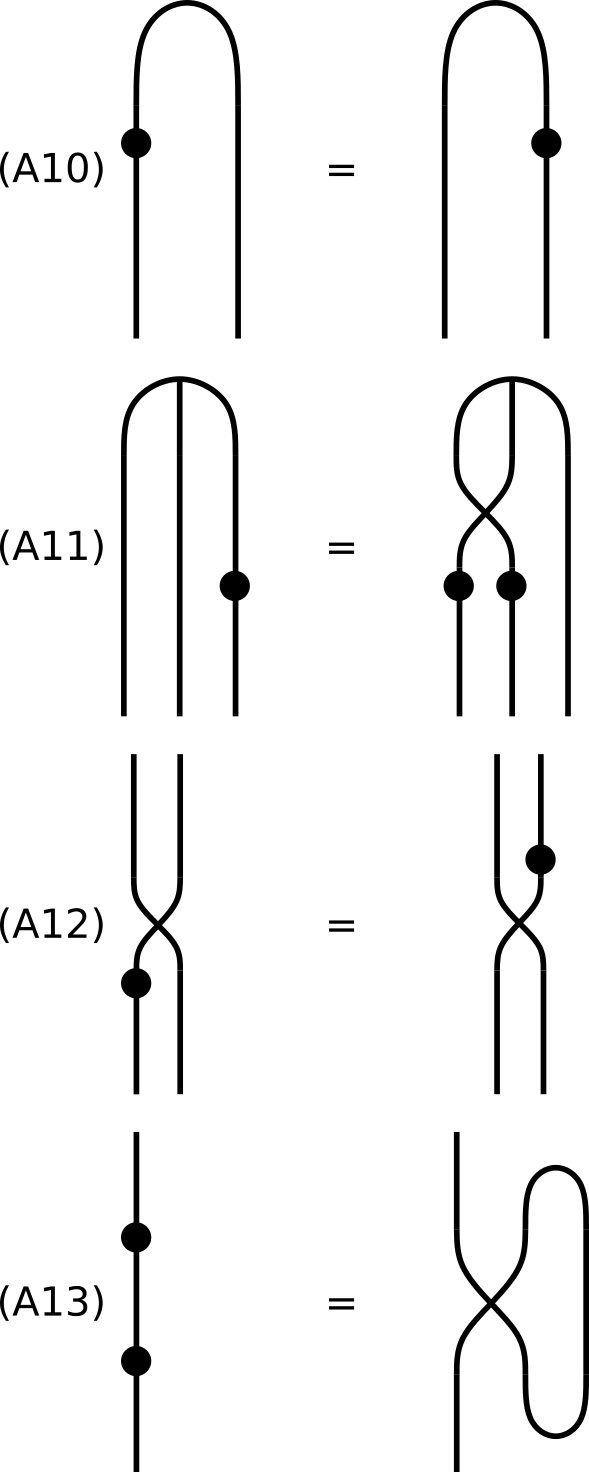}
\end{subfigure}
\caption{Ribbon diagrams for the conditions \eqref{a1}--\eqref{a13}, each due to regular homotopy or Pachner moves.}
\label{moves}
\end{figure}

\begin{figure}

	\centering
	\includegraphics[width=13cm]{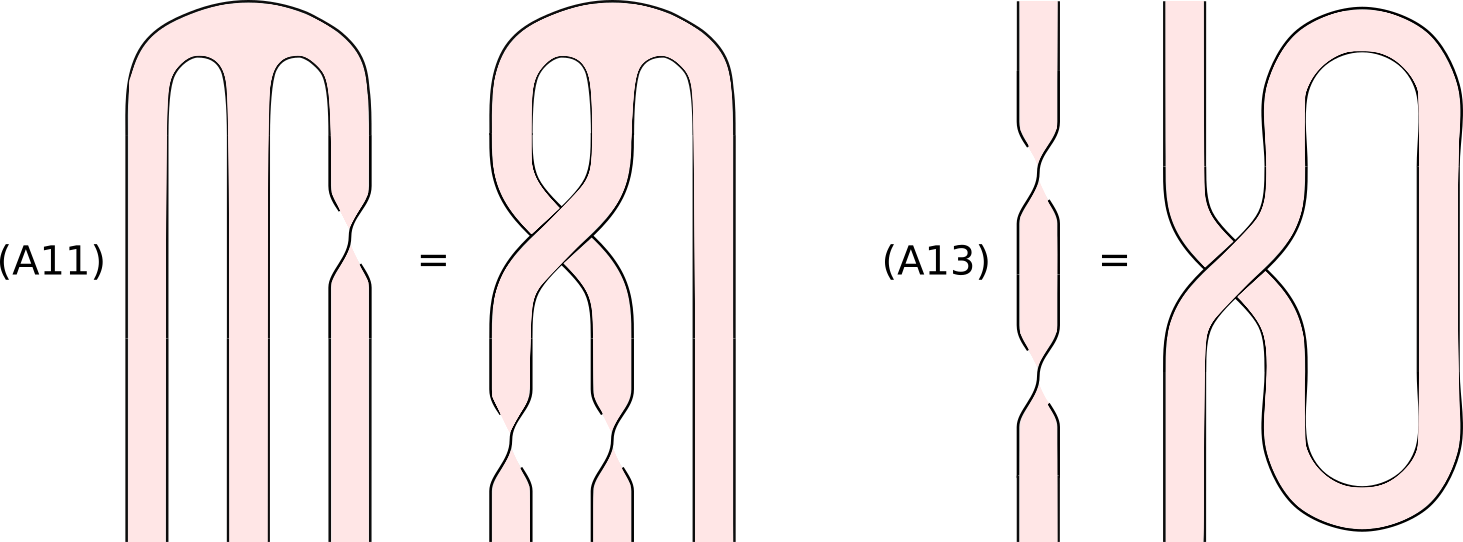}

  \caption{The moves \eqref{a11} and \eqref{a13} as ribbon diagrams.}
  \label{ribbonaxioms}
\end{figure}

If two ribbon diagrams are built from the same regular homotopy class of ribbon graphs, they are related by the set of moves depicted in Figure \ref{moves}. Moves \eqref{a1} and \eqref{a7}--\eqref{a9} are the ribbon Reidemeister moves\footnote{Note that the first ribbon Reidemeister move \eqref{a7} is weaker than the first of the usual Reidemeister moves for knots, which does not preserve the ribbon structure.} \cite{K90,FY89}. Moves \eqref{a2} and \eqref{a5}--\eqref{a6} are additional moves for graphs with nodes \cite{Y89,K05,RT90}. The moves \eqref{a10}--\eqref{a13} involve half twists and have been studied in ref. \cite{V07}. The moves\footnote{Two half twists is a full twist, and the ribbon Reidemeister moves show that a pair of full twists can be undone.} show that a left handed twist is related by regular homotopy to a sequence of three right handed twists. This means, by replacing each left handed half twist by three right handed half twists, one obtains a ribbon diagram where the half twists are all right handed. In the following, we simplify the algebra by assuming that all half twists are right handed. Two of the moves, which may be more difficult to visualize, are depicted in ribbon form in Figure \ref{ribbonaxioms}.

Any two triangulations on $M$ are related by the 2-2 \eqref{a3} and 3-1 \eqref{a4} Pachner moves \cite{P91,L99,CL14}, also depicted in Figure \ref{moves}.

\subsection{Algebraic structure}

We now show how to evaluate a partition function for a regular homotopy class of immersed surfaces. Begin with a ribbon diagram, decomposed into the five building blocks. Color the diagram by labeling the legs of each block by elements in a finite set $\mathcal{I}$. The blocks are assigned the following $\CC$-valued weights:
\begin{enumerate}
	\item Nodes labeled left to right by $a,b,c\in\mathcal{I}$ receive a weight $C_{abc}$.
    \item Caps labeled left to right by $a,b\in\mathcal{I}$ receive a weight $B_{ab}$, while cups receive a weight $B^{ab}$.
    \item Crossings labeled as in Figure \ref{blocks} by $a,b,c,d\in\mathcal{I}$ receive a weight $\lambda_{ab}{}^{cd}$.
    \item (Right handed) half twists labeled bottom to top by $a,b\in\mathcal{I}$ receive a weight $\tau_a{}^b$.
    \item Vertices\footnote{Surfaces with boundary are discussed in section \protect\ref{statespaces}. In this more general case, only \emph{internal} vertices receive a weight $R$.} of the triangulation receive a weight $R$.
\end{enumerate}

The weight of the colored diagram is the product of the weights of the pieces in its decomposition, and the partition function for a diagram is a sum of the weights of its colorings.

For the partition function to be independent of the discretization, it must be invariant under the moves of Figure \ref{moves}. By evaluating them according to our procedure, we find the following algebraic conditions:

\begin{AxAlign}
   &\text{(Snake)}&&\text{define }\eta &&&B_{ac}B^{cb}=\delta_a^b\label{a1}\\
   &\text{(Cyclicity)}&&\text{define }m &&&C_{abd}B^{dc}=B^{cd}C_{dab}\label{a2}\\
   &\text{(Pachner 2-2)}&&m\text{ associative} &&&C_{abe}B^{ef}C_{fcd}=C_{bce}B^{ef}C_{afd}\label{a3}\\
    &\text{(Pachner 3-1)}&&\eta\text{ special} &&&C_{abc}=R\,C_{ade}B^{df}C_{fbg}B^{gh}C_{ihc}B^{ei}\label{a4}\\
    \cline{1-2}
  &\text{(Crossing at a critical point)}&& &&&B_{ae}\lambda_{bc}{}^{ed}=\lambda_{ab}{}^{de}B_{ec}\label{a5}\\
  &\text{(Crossing at a node)}&& &&&\lambda_{ab}{}^{ef}C_{fcd}=C_{aeg}\lambda_{bc}{}^{ef}\lambda_{fd}{}^{ge}\label{a6}\\
  &\text{(Modified Reidemeister I)}&& &&&B^{cd}B_{ce}\lambda_{da}{}^{eb}=\lambda_{ac}{}^{bd}B^{ce}B_{de}\label{a7}\\
  &\text{(Reidemeister II)}&&&&&\lambda_{ab}{}^{ef}\lambda_{ef}{}^{cd}=\delta_a^c\delta_b^d\label{a8}\\
  &\text{(Reidemeister III)}&& &&&\lambda_{ag}{}^{di}\lambda_{bc}{}^{gh}\lambda_{ih}{}^{ef}=\lambda_{ab}{}^{gh}\lambda_{hc}{}^{if}\lambda_{gi}{}^{de}\label{a9}\\
  	\cline{1-2}
  &\text{(Twist at a critical point)}&&\eta(\mathds{1}\otimes\tau)=\eta(\tau\otimes\mathds{1}) &&&B_{ac}\tau_b{}^c=\tau_a{}^cB_{cb}\label{a10}\\
  &\text{(Twist at a node)}&&\tau m=m\lambda(\tau\otimes\tau) &&&C_{abd}\tau_c{}^d=\tau_a{}^d\tau_b{}^e\lambda_{de}{}^{fg}C_{fgc}\label{a11}\\
  &\text{(Twist at a crossing)}&&\lambda(\tau\otimes \mathds{1})=(\mathds{1}\otimes\tau)\lambda &&&\tau_a{}^e\lambda_{eb}{}^{cd}=\lambda_{ab}{}^{ce}\tau_e{}^d\label{a12}\\
  &\text{(Two half twists)}&&\tau^2=\phi &&&\tau_a{}^c\tau_c{}^b=\lambda_{ac}{}^{bd}B^{ce}B_{de}\quad (=\lambda_{ac}{}^{bd}\sigma_d{}^c=\phi_a{}^b)\label{a13}
\end{AxAlign}

The conditions \eqref{a1}--\eqref{a4} define a special Frobenius algebra $(A,m,\eta)$; that is, an a unital, associative algebra $(A,m)$ with a non-degenerate bilinear form $\eta$ satisfying the Frobenius condition $\eta(xy,z)=\eta(x,yz)$, $x,y,z\in A$, and the specialness\footnote{Ref. \protect\cite{LP06} discusses a generalization of the oriented state sum construction to non-special Frobenius algebras, where \emph{window elements} $a^{-1}\nsim 1$ are attached to vertices. In their language, we always take $a^{-1}=R\,1$ with $R\in\CC$.} condition $m\circ\eta^{-1}=R^{-1}\,1$. This algebra is defined on the vector space with basis $\{e_a\}$, $a\in\mathcal{I}$, has product $m(e_a\otimes e_b)=C_{ab}{}^ce_c$ given by associative structure coefficients $C_{ab}{}^c=C_{abd}B^{dc}$, unit $1=B^{ab}C_{bcd}B^{cd}e_a$, and non-degenerate bilinear form $\eta(e_a,e_b)=B_{ab}$. Ref. \cite{BT14} shows that the conditions \eqref{a1}--\eqref{a4} enforce the axioms of a special Frobenius algebra and, conversely, that a special Frobenius algebra defines tensors $C_{abc}$ and $B_{ab}$ that satisfy these conditions. If $\eta$ is taken to be the unique (up to $R$) \emph{symmetric} special Frobenius form, this result reduces to the familiar case studied by Fukuma, Hosono, and Kawai \cite{FHK94}.

The conditions \eqref{a5}--\eqref{a9} imply other relations like $B^{be}\lambda_{ea}{}^{cd}=\lambda_{ae}{}^{bc}B^{ed}$. The existence of a \emph{symmetric structure} $\lambda:A\otimes A\ra A\otimes A$, satisfying the axioms, is also a constraint on $\eta$. The \emph{Nakayama automorphism}
\begin{equation}\label{nakayama}
    \sigma_a{}^b=B_{ac}B^{bc},\qquad \eta(a,b)=\eta(\sigma(b),a)
\end{equation}
measures the failure of $\eta$ to be symmetric. Ref. \cite{BT14} demonstrates that conditions \eqref{a1}--\eqref{a9} imply
\begin{equation}
    B_{ac}B^{bc}=B_{ca}B^{cb},\qquad \sigma^2=1,
\end{equation}
equivalently, that $\eta$ decomposes as a sum of symmetric and antisymmetric parts. Define the \emph{full twist}
\begin{equation}\label{fulltwist}
    \phi_a{}^b=\lambda_{ac}{}^{bd}B^{ce}B_{de}=\lambda_{ac}{}^{bd}\sigma_d{}^c.
\end{equation}
Ref. \cite{BT14} argues from these conditions that $\phi$ is a Frobenius algebra automorphism; that is,
\begin{equation}\label{gradedalg}
    \phi\circ m(a\otimes b)=m(\phi(a)\otimes\phi(b)),\qquad\eta(\phi(a),\phi(b))=\eta(a,b).
\end{equation}
Moreover, $\phi$ is an involution and so defines a $\ZZ/2$-grading on $A$: on homogeneous elements,
\begin{equation}\label{grading}
    \phi(a)=(-1)^{|a|}a,\qquad |a|\in\{0,1\}.
\end{equation}
The data $(C,B,\lambda)$ satisfying these axioms is what ref. \cite{BT14} use to define their spin state sums.

Other relations like $B^{cb}\tau_c{}^a=B^{ad}\tau_d{}^b$ and $\tau_b{}^e\lambda_{ae}{}^{cd}=\lambda_{ab}{}^{fd}\tau_f{}^c$ follow from the conditions \eqref{a10}--\eqref{a13}. We will refer to the data $(C,B,\lambda,\tau)$ as a \emph{half twist algebra}. It is the input for our state sum construction.

\subsection{State spaces and bordisms}\label{statespaces}

The construction has so far focused on closed surfaces. In order to define a TQFT, it must also assign state spaces $\cA_0, \cA_1$ to one dimensional closed pin manifolds $S_0,S_1$ and linear maps $Z(M):\cA_0\ra\cA_1$ to the pin bordisms $M$ between them. Given an immersion of $M$, set up according to section \ref{pincob}, form its ribbon diagram as usual. Suppose there are $n$ edges in the triangulation of $S_0$ and $m$ in that of $S_1$. Then the state sum over internal colorings defines a map $\otimes^nA\ra\otimes^mA$. This map has a clear dependence on the triangulation, as re-triangulating may change $n$ and $m$. It is also non-invariant under regular homotopy, as crossing the external legs over each other introduces single factors of the crossing map $\lambda$. The following discussion shows that both of these problems are solved by composing each end with a certain projector.

\begin{figure}

	\centering
	\includegraphics[width=12cm]{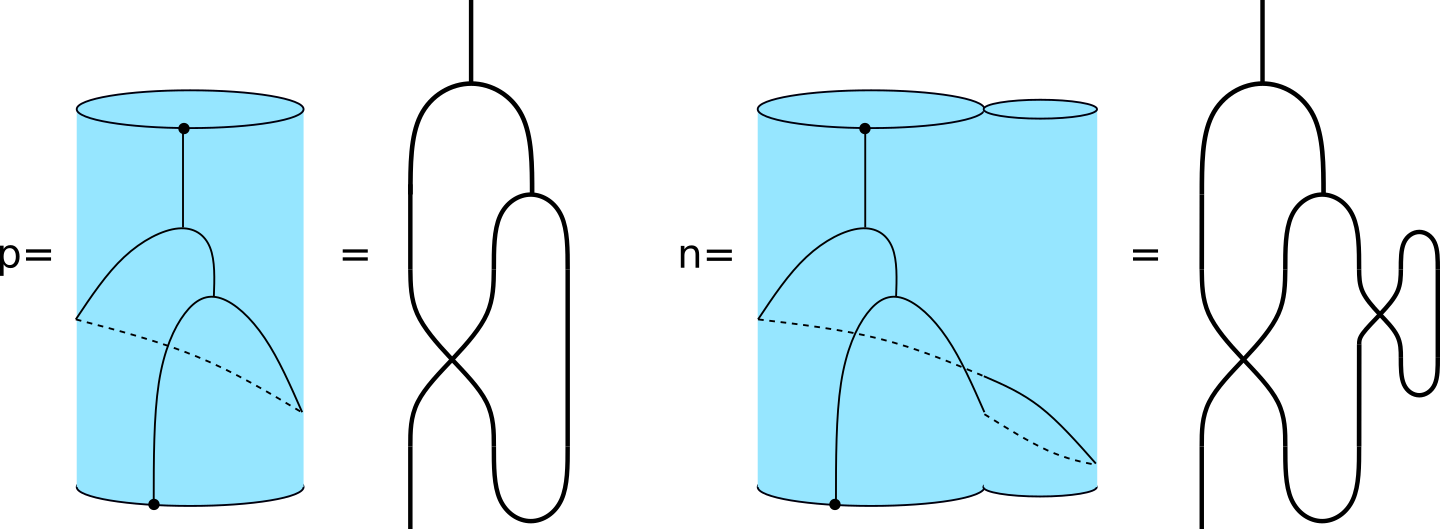}

  \captionof{figure}{Ribbon diagrams for the projectors $p$ and $n$. We have been careful to account for the half twists that appear when the ribbon turns a ``corner,'' setting them up to cancel.}
  \label{cylinder}
\end{figure}

Consider the ribbon diagrams depicted in Figure \ref{cylinder}, which arise from immersions of cylindrical topologies. One diagram corresponds to a cylinder with boundary circles of NS type, the other R.\footnote{The ribbon diagrams for cylinders of circles with rotation numbers $n,n+2$ are related by the ribbon Reidemeister moves.} Since the cylinder defines a regular homotopy between the input and output circles, they are immersed in the same way.

It has been argued by ref. \cite{BT14} (see also \cite{MS06}) that these diagrams define projectors $p$ and $n$ onto subspaces \begin{align}
&\text{im }p=\cA_{NS}=\{a\in A:m(b\otimes a)=m\circ\lambda(b\otimes a),\forall b\in A\}\label{ANS}\\
&\text{im }n=\cA_{R}=\{a\in A:m(b\otimes a)=m\circ\lambda(\phi(b)\otimes a),\forall b\in A\}\label{AR}
\end{align}
The maps assigned to other ribbon diagrams with cylindrical topology are related to these by composition with some power of $\tau$, and we will not consider them here. By gluing a copy of $p$ into each NS-type connected component of $S_0,S_1$ and a copy of $n$ into each R-type component, the map $\otimes^n A\ra\otimes^m A$ becomes
\begin{equation}
    \cZ(M):\cZ(S_0)\ra\cZ(S_1),
\end{equation}
where $\cZ(S_0)$ consists of a copy of $\cA_{NS}, \cA_R$ for each NS-type component and R-type component, respectively, and likewise for $\cZ(S_1)$. This solves the problem of triangulation-dependence.

\begin{figure}

	\centering
	\includegraphics[width=12cm]{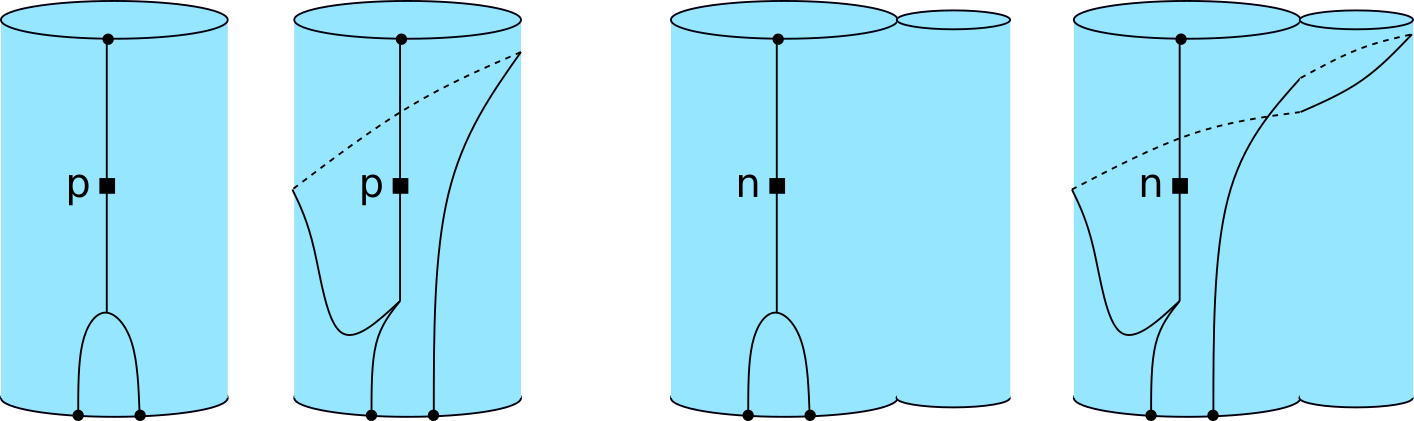}

  \captionof{figure}{Gluing independence. Since $p$ and $n$ project onto certain twisted centers of $A$, according to \eqref{ANS} and \eqref{AR}, an external leg may be pulled around the circle without affecting the state sum.}
  \label{gluingindep}
\end{figure}

One must check whether composition with $p$ and $n$ is independent of the way in which the cylindrical ribbon diagrams are glued into the cobordism. Regular homotopy has been used to push the legs of the cylindrical ribbon diagrams to the ``front'' (positive $z$-coordinate) side of the cylinders, so it must also be checked that our construction of $\cZ(M)$ is independent of the way in which this was done. Both of these checks follow from \eqref{ANS} and \eqref{AR}, which show that $p$ and $n$ are unchanged by cyclic permutation of the legs, as in Figure \ref{gluingindep}. The only ambiguity that remains is due to reordering the boundary components, which introduces factors of $\lambda$. These terms reflect the fact that the product assigned to the pair-of-pants cobordism is not commutative, but twisted-commutative. To obtain a definite $\cZ(M)$, one must fix an ordering of the boundary components; this is a characteristic of the continuum pin TQFT and not a relic of the state sum construction. For the special class of theories discussed in section \ref{supalgabk}, the product is graded-commutative with respect to the supervector structure on $\cA_{NS},\cA_R$.  In this case, $\cZ(M)$ may be interpreted as a map $\wedge_i\cZ(S^1_{0,i})\ra \wedge_i\cZ(S^1_{1,i})$ of exterior algebras, where $S^1_{0,i},S^1_{1,i}$ denote boundary components.

An axiom of (pin) TQFT requires that gluing two bordisms $M_1, M_2$ along their cut boundaries amounts to composing the linear maps assigned to them. This is true of the present construction. To see this, start by leaving off the projectors $p$ and $n$, so that the bordisms --- for some fixed discretizations --- are assigned matrices $z(M_1):\otimes^nA\ra\otimes^mA$ and $z(M_2):\otimes^mA\ra\otimes^lA$. The amplitude for the composite bordism is a sum over colorings of the internal edges of $M_1, M_2$ as well as the edges of the glued boundary, weighted the product of the weights for $M_1, M_2$. This is matrix multiplication, so $z(M_2\circ M_1)=z(M_2)z(M_1)$. To complete the argument, add back the projectors $P=p\cdots p\otimes n\cdots n$. By re-triangulation invariance, the insertion of the projectors at the glued boundary must have no effect on the state sum, so
\begin{equation}\cZ(M_2)\cZ(M_1)=Pz(M_2)PPz(M_1)P=Pz(M_2)z(M_1)P=Pz(M_2\circ M_1)P=\cZ(M_2\circ M_1).
\end{equation}

A \emph{Hermitian structure} on a pin TQFT is a sesquilinear form $\langle\cdot,\cdot\rangle$ on $\cZ(S)$ for each closed one dimensional pin manifold $S$, with respect to which $\cZ(M)$ and $\cZ(-M)$ are adjoint for any cobordism $(M,S_0,S_1)$ \cite{T10,FH16}. Here, $-M$ denotes the ``opposite'' pin cobordism from $S_1$ to $S_0$. In terms of immersed surfaces, $-M$ is obtained from $M$ by reflecting over an $xz$-plane. A \emph{unitary structure} is a Hermitian structure for which the sesquilinear form is positive definite (an inner product).

\section{Real superalgebras and the Arf-Brown-Kervaire TQFT}\label{supalgabk}

The remainder of this paper focuses on a special class of half twist algebras closely related to separable real superalgebras, the state sum models associated to which constitute a broad class of interesting examples such as the Arf-Brown-Kervaire theory. To be precise, these state sums take as a input a symmetric special Frobenius real superaglebra or, equivalently, a separable real superalgebra with a continuous parameter $\alpha$.\footnote{Sometimes we neglect $\alpha$ and speak only of the superalgebra; this is because $\alpha$'s contribution is just an Euler term.}

\subsection{Real superalgebras}\label{realsuperalgebras}

A real superalgebra is an algebra $(A_r,m)$ over $\RR$ with a linear involution $\phi:a\mapsto (-1)^{|a|}a$, with respect to which the product $m$ is equivariant, as in \eqref{gradedalg}. Superalgebras inherit the natural symmetric structure
\begin{equation}\label{svect}
    \lambda:a\otimes b\mapsto (-1)^{|a||b|}b\otimes a
\end{equation}
from the symmetric monoidal category of supervector spaces sVect. Separability\footnote{We are conflating separability and strong separability, which are equivalent conditions over $\RR$ or $\CC$.} means there is a symmetric\footnote{Here we mean ``symmetric'' in the usual sense, as a Frobenius algebra object in the symmetric monoidal category of vector spaces Vect, not that of supervector spaces sVect.} special Frobenius inner product $\eta$, unique up the nonzero real scalar $\alpha$, given by the \emph{trace form}
\begin{equation}\label{traceform}
    \eta(x,y)=\alpha\Tr[L(x)L(y)],
\end{equation}
where $L:A\ra\End(A)$ denotes left multiplication. The real algebra $A_r$ is equivalent to its complexification $A=A_r\otimes_\RR\CC$ together with an antilinear automorphism $T$ of $A$, called a \emph{real structure}, that fixes $A_r$.

By virtue of being special Frobenius, the complex algebra $A$ is separable as a superalgebra. This means it is a direct sum of simple superalgebras (``blocks''), of which there are two types: matrix algebras $\CC(p|q)$ and odd algebras $\CC(n)\otimes\CC\ell(1)$ \cite{D99}. Each block has an involutive antilinear anti-automorphism given by conjugate transposition of $\CC(p|q)$ or the $\CC(n)$ factor.\footnote{There may exist other such maps, but our construction uses this canonical one. In any basis $\{e_{ij}\}$ where $e_{ij}e_{jk}=+\,e_{ik}$, ``conjugate transposition'' is unambiguously defined as the map $e_{ij}\mapsto e_{ji}$.} The direct sum of these is a map $*$ on $A$. Its composition with the real structure is a linear involutive anti-automorphism $t=*T$.

The structures $m$, $\eta$, $\lambda$, and $\phi$ of $A_r$ extend linearly onto $A$, where the map $t$ satisfies
\begin{equation}\label{real}
    \eta(tx,ty)=\eta(x,y),\qquad tm(x\otimes y)=m(ty\otimes tx),\qquad \lambda(t\otimes\mathds{1})=(\mathds{1}\otimes t)\lambda,\qquad t^2=\mathds{1}.
\end{equation}
These relations resemble the four half twist axioms \eqref{a10}--\eqref{a13} but are not quite the same: while $t$ is $\eta$-orthogonal, $\tau$ is $\eta$-symmetric; while $t$ is an anti-automorphism, $\tau$ is a $\lambda$-twisted-automorphism; while $t$ is an involution, $\tau$ squares to $\phi$. Outside of these differences, $A$ is much like a half twist algebra: its involution $\phi$ is determined by the symmetric structure $\lambda$ as $\phi_a{}^b=\lambda_{ac}{}^{bc}$, and it is straightforward to verify that $m$, $\eta$, and $\lambda$ are compatible in the sense that they satisfy the first nine axioms \eqref{a1}--\eqref{a9}.

To make $A$ into a genuine half twist algebra, we would like to construct a half twist $\tau$, satisfying \eqref{a10}--\eqref{a13}, out of the involutive linear anti-automorphisms $t$ (associated with $T$), satisfying \eqref{real}. If $s(x)\in\{0,1\}$ is any grading of the algebra that shares an eigenbasis with $\phi$ (such as $s=0$), we may define
\begin{equation}\label{ttau}
    \tau:x\mapsto (-1)^{s(x)}i^{|x|}t(x).
\end{equation}
It is straightforward to verify that $\tau$ squares to $\phi$ and is $\eta$-symmetric. Moreover, $t$ is a $\lambda$-twisted-automorphism:
\begin{align}
\begin{split}
m\circ\lambda(\tau(x)\otimes\tau(y))&=(-1)^{|x||y|}m(\tau(y)\otimes\tau(x))\\
&=(-1)^{s(x)+s(y)}i^{|x|+|y|-2|x||y|}m(t(y)\otimes t(x))\\
&=(-1)^{s(m(x\otimes y))}i^{|m(x\otimes y)|}t\circ m(x\otimes y)\\
&=\tau\circ m(x\otimes y).
\end{split}
\end{align}
The choice of $s$ has to do with the decomposability of the state sum and is discussed in section \ref{morita}. A half twist algebra constructed from a real superalgebra is not generic. In particular, its crossing map is given by eq. \eqref{svect} and its half twist satisfies $*\tau*=\tau^{-1}$. The symmetry of $\eta$ is not an independent condition, as the special form of $\lambda$ means that the Nakayama automorphism \eqref{nakayama} is trivial.

It is worth noting at this point that our separable superalgebras come with an sesquilinear form
\begin{equation}\label{innerproduct}
	\langle x,y\rangle=\eta(*x,y).
\end{equation}
In fact, if $\alpha$ is positive, $\langle\cdot,\cdot\rangle$ is positive definite and so defines an inner product. By \eqref{traceform} it is clear that $\eta$ vanishes if $x$ and $y$ are supported on different blocks. On an even block, $\langle M,N\rangle=\alpha\Tr[M^\dagger N]$, which is positive definite. On an odd block, $\langle M\otimes\gamma^i,N\otimes\gamma^j\rangle=\alpha\,\delta^{ij}\Tr[M^\dagger N]$, which is also positive definite.

In any theory, the circles $S^1_{NS}$ and $S^1_R$ have macaroni bordisms,\footnote{Macaroni bordisms are cylinders with two ingoing boundary components. Accounting for spin structures, there are two distinct such bordisms on $S^1_R$. Choose one. The other is related by composition with a cylinder.} whose partition functions define bilinear forms $\eta_{NS}:\cA_{NS}\otimes\cA_{NS}\ra\CC$ and $\eta_R:\cA_R\otimes\cA_R\ra\CC$. Evaluating ribbon diagrams for the macaroni bordisms gives these maps in terms of the superalgebra data: $\eta_{NS}=\eta(p,p)$ and $\eta_R=\eta(n,n)$. Inserting the map $*$, as in \eqref{innerproduct}, one may define sesquilinear forms $\langle,\rangle_{NS}=\eta_{NS}(*,)$ and $\langle,\rangle_R=\eta_R(*,)$. The form on an arbitrary closed one dimensional pin manifold $S$ is given as a tensor product of these forms.

We would like to show that state sum pin TQFTs associated with real separable superalgebras are unitary in the sense of section \ref{statespaces}. It remains to check adjointness. Due to the form of $\eta$ \eqref{traceform}, this condition reads $*\cZ(M)*=\cZ(-M)^T$. In terms of ribbon diagrams in $\RR^2$, reflection across the $y$ axis must have the effect of acting on each external leg by $*$. The conditions on each building blocks read
\begin{equation}\label{adjoint}
    *m(*a\otimes*b)=m(b\otimes a),\quad\eta(*a\otimes *b)=\eta(b\otimes a),\quad(*\otimes *)\lambda(*a\otimes *b)=\lambda(b\otimes a),\quad *\tau(*a)=\tau^{-1}(a).
\end{equation}
The first condition follows from the fact that $*$ is an anti-automorphism, the second and third from symmetry of $\eta$ \eqref{traceform} and $\lambda$ \eqref{svect}, and the fourth from the antilinearity of $*$ and the $i$ factor in \eqref{ttau}. Unitarity also requires $R\in\RR$, which follows from $\alpha\in\RR$. Therefore theories associated to real separable superalgebras are unitary.

\begin{figure}

	\centering
	\includegraphics[width=13cm]{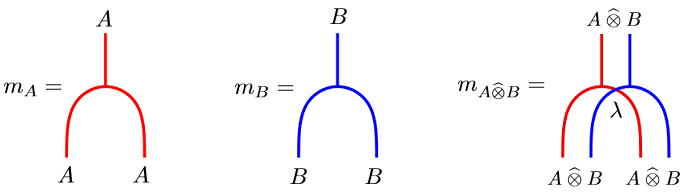}

  \captionof{figure}{A diagrammatic representation of the supertensor product of superalgebras $A$ and $B$.}
  \label{AB}
\end{figure}

A useful construction on superalgebras $A,B$ is the \emph{supertensor product} $A\fotimes B$. This superalgebra has underlying vector space $A\otimes B$ with grading $\phi_{A\fotimes B}=\phi_A\fotimes\phi_B$ and associative product
\begin{equation}\label{supertensor}
    (a\fotimes b)(a'\fotimes b')=(-1)^{ba'}aa'\fotimes bb',\qquad\text{ i.e. }\quad m_{A\fotimes B}=(m_A\fotimes m_B)(1\fotimes\lambda_{BA}\fotimes 1),
\end{equation}
where $\lambda_{AB}:A\fotimes B\ra B\fotimes A$ is the symmetric structure of sVect \eqref{svect}. The special symmetric Frobenius form is $\eta_{A\fotimes B}=(\eta_A\fotimes\eta_B)(1\fotimes\lambda\fotimes 1)$. It is helpful to interpret the product rule \eqref{supertensor} diagrammatically. In Figure \ref{AB}, the products on $A$ and $B$ are represented by trivalent nodes of red and blue lines, respectively. The product on $A\fotimes B$ has a red-blue crossing, contributing the sign $\lambda$. More generally, one may consider diagrams that consist of a red ribbon diagram superimposed on a blue ribbon diagram such that the usual regularity conditions are met. Color the red diagram by basis elements $e_a$ of $A$ and the blue diagram by basis elements $f_i$ of $B$. The weight of this double coloring is the weight of the red coloring, according to $A$, times the weight of the blue coloring, according to $B$, times signs $|e_a||f_i|$ at each red-blue crossing. It is invariant under the usual moves \eqref{a1}-\eqref{a13} of each of the red and blue diagrams. Due to the graded products on $A$ and $B$, the weight is also invariant under these same moves  where some of the ribbons are red and some are blue. In particular, the weight is unchanged by pulling a red-blue crossing across a critical point, node, or half twist, and satisfies colored versions of the ribbon Reidemeister moves. This sort of representation will prove useful in section \ref{morita} when we discuss the state sum for $A\fotimes B$. It is worth mentioning that the supertensor product is the monoidal product of superalgebras when they are regarded as algebra objects in sVect. In this language, the colored moves are related to the graphical calculus of symmetric monoidal categories.

\subsection{Example: Clifford algebras}\label{clifford}

In this section, we define the Clifford algebras $C\ell_{p,q}\RR$ and $C\ell_n\CC$ and discuss their associated half twist algebras, from which one can extract the state sum data $(C,B,\lambda,\tau)$. As will be shown in section \ref{morita}, the significance of these examples is that they generate all theories associated to separable real superalgebras.\footnote{We leave open the question of whether there exist pin-TQFTs that do not arise via our state sum construction.}

The real Clifford algebra $A=C\ell_{p,q}\RR$ is generated by anticommuting elements $\gamma_1,\ldots,\gamma_p$ with $\gamma_j^2=+1$ and $\gamma_{p+1},\ldots,\gamma_{p+q}$ with $\gamma_j^2=-1$. It has a basis $\{\gamma_1^{N_1}\cdots\gamma_n^{N_n}\}$ for $N_j=0,1$, $n=p+q$. The form $\eta=\epsilon\circ m$ given by the counit
\begin{equation}
    \epsilon(\gamma_1^{N_1}\cdots\gamma_n^{N_n})=\left\{\begin{array}{lr}\alpha\,2^{n/2}&N_j=0, \forall j\\0&\text{else}\end{array}\right.
\end{equation}
is Frobenius, symmetric, and special with $R=\alpha\,2^{-n/2}$. The grading is given by the standard involution
\begin{equation}\label{cliffordinvolution}
    \phi(\gamma_1^{N_1}\cdots\gamma_n^{N_n})=(-1)^{\sum_jN_j}\gamma_1^{N_1}\cdots\gamma_n^{N_n}.
\end{equation}
For the element $x=\gamma_1^{N_1}\cdots\gamma_n^{N_n}$, let $\{x\}=\sum_j N_j$, which is to say $|x|=\{x\}\text{ mod }2$.

The corresponding half twist algebra is defined on the complexification $C\ell_{p+q}\CC=C\ell_{p,q}\RR\otimes_\RR\CC$, which comes with a real structure $T$ that fixes the $\gamma$-basis and complex conjugates its coefficients. Let us define new generators $\Gamma_j=\gamma_j$ for $1<j\le p$ and $\Gamma_j=i\gamma_j$ for $p<j\le p+q$, so that $\Gamma_j^2=+1$. The basis element $x=\Gamma_1^{N_1}\cdots\Gamma_n^{N_n}$ has $T$-eigenvalue $(-1)^{|x|_q}$, where $|x|_q=\sum_{i>p}N_i\text{ mod }2$. It remains to construct the half twist $\tau$. The Clifford algebra has a natural Hermitian structure given by the conjugate transpose map
\begin{equation}
    *(\Gamma_1^{N_1}\cdots\Gamma_n^{N_n})=\Gamma_n^{N_n}\cdots\Gamma_1^{N_1}=(-1)^{\{x\}(\{x\}-1)/2}\Gamma_1^{N_1}\cdots\Gamma_n^{N_n}.
\end{equation}
The composition $t=*T$ fails the condition \eqref{a13}; however, it can be corrected, as in eq. \eqref{ttau} with $s=0$. Define
\begin{equation}
    \tau(x)=i^{|x|}t(x)=i^{|x|}(-1)^{|x|_q}(-1)^{\{x\}(\{x\}-1)/2}x=i^{\{x\}}(-1)^{|x|_q}x.
\end{equation}
The general discussion in section \ref{realsuperalgebras} shows that the half twist axioms are satisfied.

The complex Clifford algebra $C\ell_n\CC$ also appears as a real superalgebra generated by anticommuting elements $\Gamma_1,\cdots,\Gamma_n$ with $\Gamma_j^2=+1$ and central element $\imath$ with $\imath^2=-1$.\footnote{This algebra is graded-isomorphic to one with $\tilde\Gamma_j^2=-1$ for some $j$ by the identification $\tilde\Gamma_j=\Gamma_j\imath$.} On basis elements $\Gamma_1^{N_1}\cdots\Gamma_n^{N_n}\imath^\cM$, the counit is $\alpha\,2^{(n+2)/2}$ if $N_j=\cM=0$ and $0$ otherwise. The form $\eta=\epsilon\circ m$ is Frobenius, symmetric, and special with $R=\alpha\,2^{-n/2}$. The central element $\imath$ is $\phi$-even, while the $\Gamma_j$ are $\phi$-odd, so $|x|=\{x\}\text{ mod }2$ where $\{x\}=\sum_jN_j$. The complexification $\CC\ell_n\CC\otimes_\RR\CC$ has real structure $T$ that fixes the $\Gamma_j$ and $\imath$. The structure $*$ is again given by conjugate transposition. According to \eqref{ttau} with $s(x)=\cM$, the half twist is the composition $\tau(x)=(-1)^\cM i^{\{x\}}x$.

\subsection{State sum for the Arf-Brown-Kervaire TQFT}\label{abkss}

The pin state sum construction discussed in section \ref{ribht} amounts to choosing a discretization of a pin surface $M$, building an associated ribbon diagram, and performing a weighted sum over colorings of the ribbon diagram. While this construction bears some resemblance to the state sums of Novak and Runkel \cite{NR15}, our approach to discretizing the (s)pin structure -- based on immersions rather than markings -- introduces a crucial difference: the existence of crossing elements means that a coloring of the ribbon diagram (in the plane) is not in general realized by a coloring of the ribbon graph (in the surface) projected onto the plane. For the present purpose of computing the state sum of the Arf-Brown-Kervaire theory, this difference is an obstacle, though one that can be avoided by restricting to the special class of half twist algebras discussed earlier in this section.

\begin{figure}

	\centering
	\includegraphics[width=15cm]{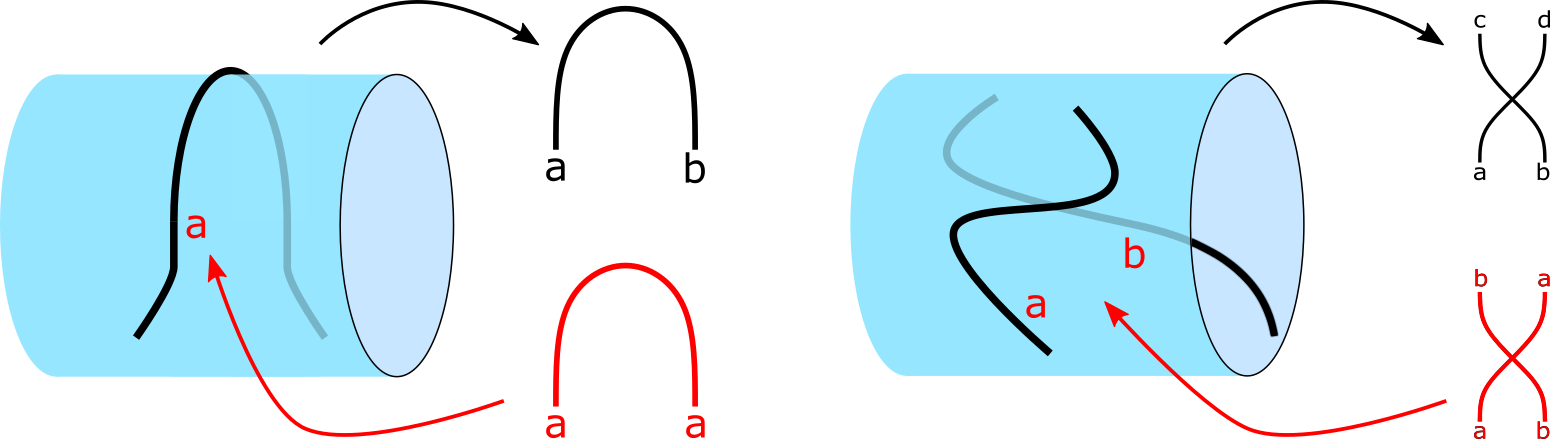}

  \captionof{figure}{If $B$ is symmetric, any coloring of the cap that has nonzero amplitude arises from a coloring of an edge in the ribbon graph. Likewise, if $\lambda$ is of the form $\lambda:a\otimes b\mapsto\lambda(a,b)b\otimes a$, any coloring of the crossing that has nonzero amplitude arises from a coloring of two edges in the ribbon graph.}
  \label{graphcolor}
\end{figure}

The state sum associated to a separable real superalgebra has the special property that it can be written as a sum over colorings \emph{of the graph} dual to the triangulation of $M$. These colorings are a special type of coloring of the ribbon diagram where all segments of a ribbon from node to node have the same label, as in Figure \ref{graphcolor}. A pin state sum localizes to these colorings if the amplitudes for all other colorings vanish. This means that $B$ is symmetric and there is a basis of $\tau$ eigenstates in which $\lambda_{ab}{}^{cd}=\lambda(a,b)\delta_a^d\delta_b^c$ for some values $\lambda(a,b)\in\CC$. By \eqref{a5} and \eqref{a8}, $\lambda(a,b)=\lambda(b,a)\in\{\pm 1\}$, and by definition of the full twist $\lambda(a,a)=(-1)^{|a|}$. The half twist algebra associated to separable real superalgebra satisfies these conditions with $\lambda(a,b)=(-1)^{|a||b|}$. The collection of edges labeled by $\phi$-odd basis elements  forms a $1$-chain $x$ with $\ZZ/2$ coefficients for the triangulation of $M$. Since the product $m$ is $\phi$-equivariant \eqref{gradedalg}, a coloring contributes zero amplitude to the state sum unless the number of odd labels surrounding each node of the graph is even; that is, unless $x$ is a cycle. Thus the sum over colorings reduces to a sum over cycles $x$:
\begin{equation}
    \cZ=\sum_{x\in Z_1(M;\ZZ/2)}\cZ(x).
\end{equation}

Consider the half twist algebra $A$ corresponding to $C\ell_{1,0}\RR$. It is spanned by $1$ and the $\phi$-odd generator $\Gamma$ with $\Gamma^2=+1$. In this basis, the tensor $B_{ab}$ is $\alpha\sqrt{2}\,\delta_{ab}$, while $C_{abc}=C_{ab}{}^dB_{dc}$ is $\alpha\sqrt{2}$ if an $|a|+|b|+|c|=0\mod 2$ and $0$ otherwise. The half twist has $\tau(1)=1$, $\tau(\Gamma)=i\Gamma$. The constant $R$ is $\alpha/\sqrt{2}$.

Each cycle $x$ is represented by a collection $\{\gamma_i\}_i$ of disjoint loops in the graph. Let us first consider the case of a single loop $\gamma$. Form a ribbon diagram and assign a weight to $\gamma$ using the data of the half twist algebra. Without loss of generality, take the legs of each $C$ to point downward and those of each $B$ upward. The tensors $C_{abc}$ and $B^{ab}$ contribute $\alpha\sqrt{2}$ and $(\alpha\sqrt{2})^{-1}$, respectively, since there are an even number of $\Gamma$ labels at each node, cap, and cup. Since the number of $C$'s is the number $|V|$ of vertices of the graph and the number of $B$'s is the number $|E|$ of edges, these contributions give an overall factor of $(\alpha\sqrt{2})^{|V|-|E|}$. Each half twist traversed by $\gamma$ contributes $i$, while each self-crossing of $\gamma$ contributes $\lambda_{\Gamma\Gamma}{}^{\Gamma\Gamma}=-1$. Therefore, the contribution of $\gamma$ to the state sum is $i^{\tilde q(\gamma)}$, where $\tilde q$ counts the number of half twists plus twice the number of crossings. It was observed in section \ref{quadsubsec} eq. \eqref{quadlocal} that this $\tilde q$ is the quadratic enhancement associated to the pin structure on $M$. Now allow for multiple loops. If the images of distinct loops intersect, they must do so at an even number of points, so the factor due to their crossing vanishes. The contribution to the state sum is $i^{\sum_j\tilde q(\gamma_j)}$. Since the loops are disjoint and so have intersection number zero, it follows from \eqref{quad} that the exponent is $\sum_j\tilde q(\gamma_j)=q(x)$, the quadratic enhancement evaluated on the cycle $x$ associated to $\{\gamma_j\}_j$. The contributions of two homologous chains differ by that of a boundary, which must be $i^{q(x)}=1$. This means that the sum over $x$ reduces to a sum over homology classes $[x]$ times the number of boundaries. This number is $2^{|F|-1}$ where $|F|$ is the number of faces of the graph.\footnote{Assuming $M$ is connected, the boundary map on $2$-cells has a two element kernel.} The full state sum is
\begin{align}\begin{split}\label{cl1abk}
\cZ_{C\ell_{1,0}\RR}(M,s)&=\left(\alpha/\sqrt{2}\right)^{|F|}\left(\alpha\sqrt{2}\right)^{|V|-|E|}\sum_{x\in Z_1(M;\ZZ/2)}e^{i\pi q_s(x)/2}\\
&=\frac{\alpha^{\chi(M)}}{\sqrt{2^{2-\chi(M)}}}\sum_{[x]\in H_1(M;\ZZ/2)}e^{i\pi q_s([x])/2}\\
&=\alpha^{\chi(M)}\text{ABK}(M,s),
\end{split}\end{align}
since $|V|-|E|+|F|$ is the Euler characteristic $\chi(M)$ and $2^{2-\chi(M)}=|H_1(M;\ZZ/2)|$.

Using the expressions \eqref{ANS} and \eqref{AR}, we find $\cA^{NS}_{C\ell_{1,0}\RR}=\CC^{1|0}$, spanned by $1$, while $\cA^R_{C\ell_{1,0}\RR}=\CC^{0|1}$, spanned by $\Gamma$. In other words, the NS sector is even (as always), while the R sector is odd (unlike the trivial theory).

Here is a good place to discuss the theory associated to the real superalgebra $C\ell_1\CC$. It is convenient to work in a basis of complex central idempotents $E_\pm=(1\pm i\imath)/2$ and elements $\Gamma E_\pm$. In this basis, $B_{ab}$ is $\alpha\sqrt{2}\,\delta_{ab}$, while $C_{abc}$ vanishes if the three $\pm$ indices do not agree or if there are an odd number of $\Gamma$'s and is otherwise $\alpha\sqrt{2}$. The half twist exchanges $E_+$ with $E_-$ and $\Gamma E_+$ with $\Gamma E_-$ while multiplying the latter two by $i$. This means that, if any loop in the ribbon diagram has an odd number of half twists, there is no way to color the edges such that the amplitude is nonzero. This happens if and only if $M$ is nonorientable; thus, the partition function vanishes on nonorientable surfaces. For orientable surfaces, it is always possible to remove all half twists from the ribbon diagram. Then, for colorings with nonzero amplitude, either all of the edges are labeled by $E_+,\Gamma E_+$ or they are all labeled by $E_-,\Gamma E_-$. In each case, such colorings are given by disjoint loops labeled by $\Gamma E$ with all other edges labeled by $E$. As above, these configurations contribute factors of $i^{q(x)}$. The contributions of the $B$ and $C$ tensors are the same as before. In total,
\begin{equation}
    \cZ_{C\ell_1\CC}(M,s)=\left\{\begin{array}{lr}2\,\alpha^{\chi(M)}\text{Arf}(M,s)&M\text{ orientable}\\0&M\text{ nonorientable}\end{array}\right.
\end{equation}
The factor of $2$ comes from the equal contributions of the $E_+,\Gamma E_+$ sector and the $E_-,\Gamma E_-$ sector, and $\text{Arf}(M,s)$ denotes the Arf invariant for the spin structure induced by the orientations and pin structure on $M$ \cite{A41,J80}; it is the restriction of the ABK invariant to orientable surfaces. One may compute the state spaces $\cA_{C\ell_1\CC}^NS=\CC^{2|0}$, spanned by $1$ and $\imath$, and $\cA_{C\ell_1\CC}^R=\CC^{0|2}$, spanned by $\Gamma$ and $\Gamma\imath$.

The vanishing of the partition function on nonorientable surfaces reflects the fact that the time reversal symmetry of the corresponding lattice model has been broken. This interpretation is also compatible with the two dimensional state spaces, which appear as ground state degeneracies in the lattice model.

\subsection{Decomposability, stacking, and Morita equivalence}\label{morita}

A TQFT $\cZ$ is said to be \emph{decomposable} if there exist TQFTs $\cZ_1,\cZ_2$ such that $\cZ\simeq\cZ_1\oplus\cZ_2$ on all spaces and cobordisms. The previous subsection demonstrated how the data of a separable real superalgebra $A$ defines a pin TQFT $\cZ_A$. We now argue that if $A$ decomposes as $A_1\oplus A_2$ the TQFT $\cZ_A$ decomposes as $\cZ_{A_1}\oplus\cZ_{A_2}$. This result motivates us to restrict our attention to indecomposable separable (a.k.a. simple) algebras.

It is clear that the circle state spaces, found in section \ref{statespaces} to be certain twisted centers of $A$, decompose as $\cA_{NS}=\cA_{NS,1}\oplus\cA_{NS,2}$ and $\cA_R=\cA_{R,1}\oplus\cA_{R,2}$. Thus $\cZ(S)\simeq\cZ_1(S)\oplus\cZ_2(S)$. A coloring of a ribbon diagram by elements in a basis of $A_1\oplus A_2$ has zero amplitude unless either all of the labels (internal and external) are from $A_1$ or they are all from $A_2$. This is the case because it holds for the building blocks $C$, $B$, and $\tau$. Therefore, $\cZ$ acts as $\cZ_1(M)$ on the subspaces $\cZ_1(S)$ and as $\cZ_2(M)$ on $\cZ_2(S)$, so $\cZ(M)\simeq \cZ_1(M)\oplus\cZ_2$(M), as claimed. In particular, when $M$ is a closed surface, $\cZ(M)=\cZ_1(M)+\cZ_2(M)\in\CC$.

The converse --- that indecomposability of $A$ implies that of $\cZ_A$ --- of the statement above is not generally true for $\cZ_A$ built out of $A$ with a half twist of the form of eq. \eqref{ttau}; however, it holds for the examples considered in section \ref{abkss} due to our careful choices of the grading $s$. The careful choice of $s$ for generic $A$ is the following. Decompose $A_r$ as a direct sum of Clifford algebras tensored with matrix algebras and choose $s=0$ on each real Clifford algebra, $s=\cM$ on each complex Clifford algebra, and $s=0$ on each matrix algebra. The complex algebra $A$ splits into blocks by orthogonal central idempotents $E_i$. With these choices, $\tau$ fixes an $E_i$ if and only if $T$ does.\footnote{In the example of $C\ell_1\CC$, the elements $E_\pm$ are fixed by neither $\tau$ nor $T$ when $s=M$ but are fixed by $\tau$ when $s=0$.} The meaning of $T$ fixing an $E_i$ is that $A_r$ decomposes along this block, while the meaning of $\tau$ fixing an $E_i$ is that the state sum decomposes. This is because, for colorings with nonzero weight, each of the three edges at a node must be colored in a single block, and so, unless $\tau$ exchanges blocks between nodes, the coloring of all edges of the ribbon diagram must be in a single block.

There is another operation on pin TQFTs called \emph{stacking}. The result of stacking $\cZ_1$ with $\cZ_2$ is the theory defined by the graded tensor product $\cZ\simeq\cZ_1\fotimes\cZ_2$. We now argue that $\cZ_{A\fotimes B}\simeq\cZ_A\fotimes\cZ_B$.

Recall that $\cA_{NS}=\text{im }p$ \eqref{ANS} and $\cA_R=\text{im }n$ \eqref{AR}. If $a\in\cA_{NS},b\in\cB_{NS}$, then for all $a\in A,b\in B$,
\begin{equation}
    (a\fotimes b)(a'\fotimes b')=(-1)^{ba'}aa'\fotimes bb'=(-1)^{ba'+aa'+bb'}a'a\fotimes b'b=(-1)^{(a+b)(a'+b')}(a'\fotimes b')(a\fotimes b),
\end{equation}
so $a\fotimes b\in\text{im }p_{A\fotimes B}$. The same argument shows the converse. Similarly, if $a\in\cA_R,b\in\cB_R$,
\begin{equation}
    (a\fotimes b)(a'\fotimes b')=(-1)^{(a+b)(a'+b')+(a'+b')}(a'\fotimes b')(a\fotimes b).
\end{equation}
Therefore, $\cZ_{A\fotimes B}(S^1_\alpha)\simeq\cZ_A(S^1_\alpha)\fotimes\cZ_B(S^1_\alpha)$ for $\alpha=\text{NS, R}$. On a one dimensional closed pin manifold,
\begin{equation}
    \cZ_{A\fotimes B}(S)=\bigfotimes_i\cZ_{A\fotimes B}(S^1_i)=\bigfotimes_i\cZ_A(S^1_i)\fotimes \cZ_B(S^1_i),
\end{equation}
which is isomorphic to $\cZ_A(S)\fotimes\cZ_B(S)$ by a sign arising from the rule \eqref{svect}. Therefore $\cZ_{A\fotimes B}\simeq\cZ_A\fotimes\cZ_B$ on the level of state spaces. Note that this argument demonstrates that the supertensor product, rather than the ordinary tensor product, is the correct stacking operation.

\begin{figure}

	\centering
	\includegraphics[width=13cm]{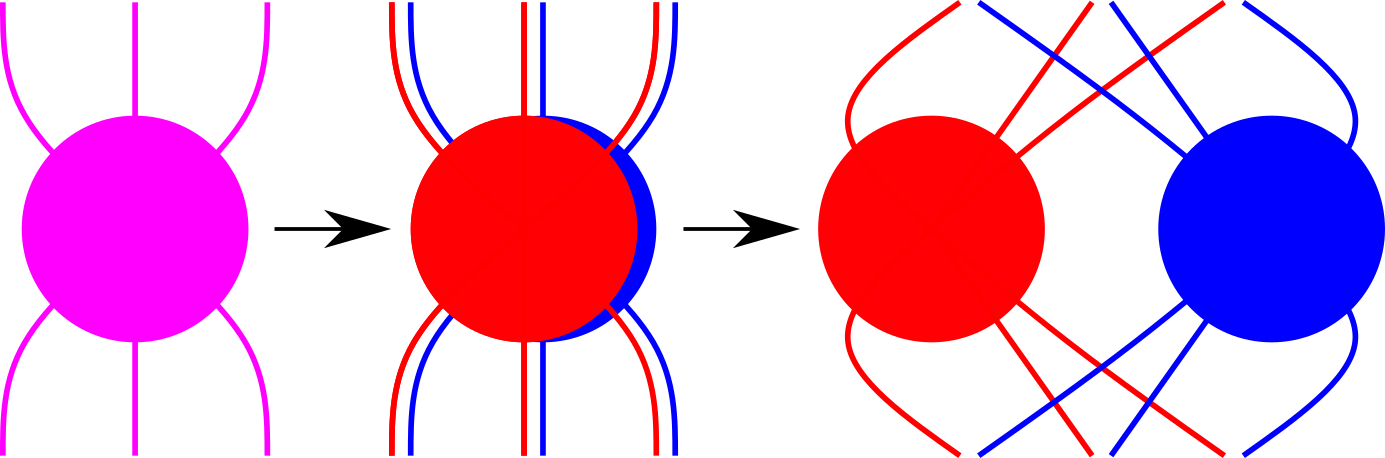}

  \captionof{figure}{A ribbon diagram for the supertensor product algebra $A\fotimes B$ (purple) may be split into a ribbon diagram for $A$ (red) superimposed on a ribbon diagram for $B$ (blue). Then they may be separated.}
  \label{stacking}
\end{figure}

The state sum for $\cZ_{A\fotimes B}$ is given by a sum over colorings of a ribbon diagram by basis elements $e_a\fotimes f_i$. One may represent these colorings as follows. Add to the ribbon diagram (in red) a copy of itself (in blue), shifted a small distance in the $x$-direction, as in Figure \ref{stacking}. The weight of this red-blue diagram, discussed in section \ref{realsuperalgebras}, reproduces the weight \eqref{supertensor} at nodes as well as the correct weights for the other building blocks in $A\fotimes B$. Now observe that the two diagrams may be pulled apart. This is allowed due to red-blue versions of the half twist axioms leaving the weight invariant. If $M$ is closed, we are done, as the weights for the $A\fotimes B$ theory are the products of those of the $A$ and $B$ theories. If $M$ has cut boundaries, we may assume that each connected component of the boundary has a single leg. Pulling apart the diagrams costs signs due to the crossings of these external legs, but these signs are precisely those in the isomorphism $\cZ_{A\fotimes B}(S)\simeq\cZ_A(S)\fotimes\cZ_B(S)$. We conclude that $\cZ_{A\fotimes B}\simeq\cZ_A\fotimes\cZ_B$ on the level of amplitudes as well.

Two superalgebras $A, B$ are said to be \emph{(graded) Morita equivalent} if their categories of graded modules are equivalent. When $A, B$ are simple, they each have a unique simple graded module (up to isomorphism, including parity change) \cite{D99}, and so Morita equivalence means that the superalgebras of module endomorphisms (the ``commutants'') of these modules are isomorphic. This relation between simple superalgebras is known as Brauer-Wall equivalence \cite{W64,D99}; another formulation, more useful for our purposes, says that $A, B$ are equivalent if they are related by stacking with matrix algebras; that is, if $A\fotimes\RR(p|q)\simeq B\fotimes\RR(p'|q')$ for some $p, q, p', q'\in\NN$ \cite{L05}. It is easy to see that the operation of stacking is compatible with this equivalence, so that one may speak of stacking equivalence classes: $[A]\fotimes[B]\simeq[A\fotimes B]$. It is worth emphasizing that the state sum construction takes as input a real superalgebra; forgetting the graded structure identifies many of these (and their Morita classes), as does complexifying and forgetting the real structure.

It will be shown in section \ref{invertible} that the pin TQFT corresponding to the algebras $\RR(p|q)$, with $\alpha=1$, is the unit in the monoid of pin TQFTs under stacking; in particular, it has state spaces $\cZ(S^1_{NS})=\cZ(S^1_R)=\CC^{1|0}$ and partition function $\cZ(M)=1$ for any closed pin surface $M$. This fact justifies the conclusion that Morita equivalent algebras $A\sim B$ define the same TQFT, $\cZ_A\simeq\cZ_B$, up to an Euler term.

The Morita-invariance of the state sum construction motivates us to focus on certain convenient representatives from each Morita class. There are ten Morita classes of simple real superalgebras \cite{W64}. Eight of them are central simple and form a group $\ZZ/8$ under stacking. The real Clifford superalgebra $C\ell_{p,q}\RR$ --- discussed in section \ref{clifford} --- lives in the class labeled by its signature $p-q\text{ mod }8$ \cite{ABS64}. The remaining two Morita classes are non-central and do not have inverses under stacking. They are represented by the complex Clifford superalgebras $C\ell_n\CC$, with $n\text{ mod }2$ being Morita invariant. In light of the result of section \ref{abkss} that the $C\ell_{1,0}\RR$ theory has partition function $\text{ABK}$, our discussion of stacking and Morita equivalence means that the algebra $C\ell_{p,q}\RR$ has partition function $\text{ABK}^{p-q}$.

\subsection{Invertible pin TQFTs}\label{invertible}

An invertible pin TQFT is one whose state spaces are one dimensional and whose partition functions on closed pin spacetimes are nonzero. Unitary\footnote{The assumption of unitarity is crucial. The non-unitary theory built from $A=\RR$ with $\alpha=-1$ has the same cobordism-invariant partition functions $(-1)^\chi=(-1)^{w_2}=(-1)^{w_1^2}=\text{ABK}^4$ as the unitary theory built from $A=C\ell_{4,0}\RR$ with $\alpha=+1$; the two theories are distinguished on the macaroni cobordism.} invertible theories\footnote{To be precise, we mean unitary invertible theories with values in supervector spaces.} have a special property \cite{Y18,FM04}: not only are they completely determined by their partition functions on closed pin manifolds, these partition functions must be a cobordism invariant --- a power of the ABK invariant --- times an Euler term $\alpha^\chi$ for $\alpha\in\RR_{>0}$.\footnote{Invertible pin TQFTs \emph{do not} generate a complete set of pin diffeomorphism invariants, as the bounding torus and bounding Klein bottle cannot be distinguished: they have both $\text{ABK}$ and $\chi$ trivial.} In particular, if $\cZ(S^2)=\alpha^2=1$, the partition functions are cobordism-invariant and multiplicative under the appropriate notion of connect sum. These theories have been constructed as extended TQFTs in ref. \cite{DG18}. Since $\text{ABK}^k(\RR P^2_1)=\exp(k\pi i/4)$ and $\text{ABK}^8=1$, the partition function on $\RR P^2_1$ (alternatively, $\RR P^2_7$) determines $k$ and therefore the full pin TQFT. In the following, we will compute the partition functions of $\RR P^2_1$ for the theories associated to the real superalgebras $\RR(p|q)$ and $C\ell_{p,q}\RR$ and find that they are $+1$ and $\exp((p-q)\pi i/4)$, respectively, up to Euler terms. Since these theories are invertible and unitary, this demonstrates that the state sum for matrix algebras is trivial --- as claimed in section \ref{morita} --- while that for $C\ell_{p,q}\RR$ is the $\text{ABK}^{p-q}$ theory --- in agreement with the findings of section \ref{abkss}.

\begin{figure}

	\centering
	\includegraphics[width=15cm]{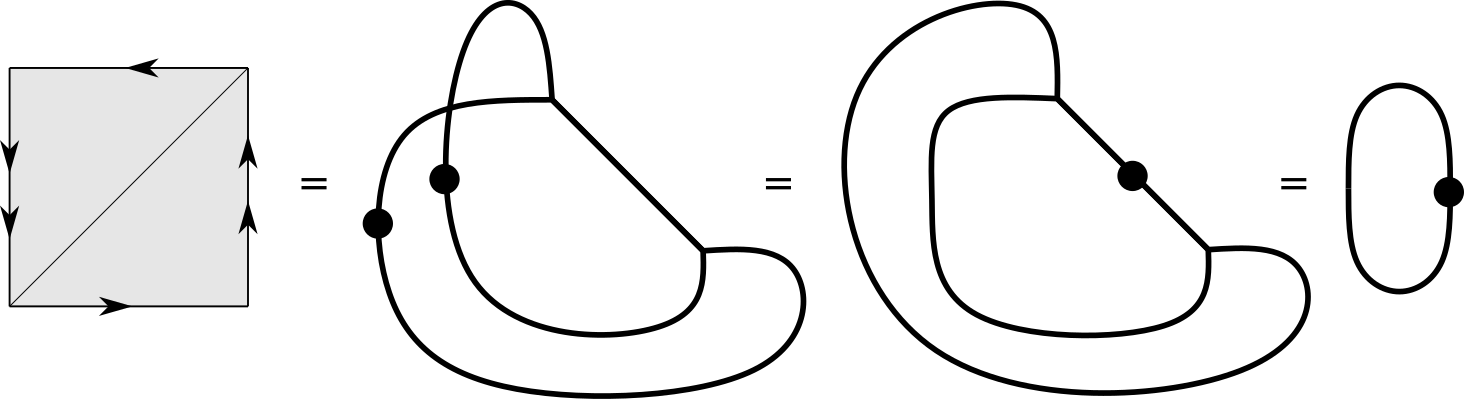}

  \captionof{figure}{A ribbon diagram for $\RR P^2$ is obtained from the graph dual to a triangulation of its fundamental square and then simplified using the moves \eqref{a11} and \eqref{a4}.}
  \label{RP2}
\end{figure}

A ribbon diagram for $\RR P^2_1$ is depicted in Figure \ref{RP2}. It evaluates to
\begin{equation}
	\cZ(\RR P^2_1)=R\,\eta(\mathds{1}\otimes\tau)\eta^{-1}.
\end{equation}
The matrix algebra $\RR(p|q)$ is spanned by a basis of matrices $e_{ij}$ with $0<i,j\le p+q=n$. The trace form is
\begin{equation}
	\eta(e_{ij},e_{kl})=\alpha\Tr[e_{ij}e_{kl}]=\alpha\,\delta_{jk}\delta_{il},\qquad \eta^{-1}=\alpha^{-1}\sum_{i,j}e_{ij}\otimes e_{ji},\qquad R=\alpha/n.
\end{equation}
Let $|i|$ be $1$ if $i>p$ and $0$ otherwise. The grading on $\RR(p|q)$ is given by $|e_{ij}|=|i|+|j|-|i||j|$. $T$ acts trivially in this basis, and $\RR(p|q)$ has a Hermitian structure given by conjugate transposition: $*e_{ij}=e_{ji}$. Therefore, by the discussion in section \ref{realsuperalgebras}, the half twist is $\tau(e_{ij})=i^{|i|+|j|+|i||j|}e_{ji}$. Then compute
\begin{equation}
	\cZ_{\RR(p|q)}(\RR P^2_1)=\frac{1}{n}\sum_{i,j}\eta(e_{ij}\otimes\tau(e_{ji}))=\frac{1}{n}\sum_{i,j}i^{|i|+|j|+|i||j|}\eta(e_{ij}\otimes e_{ij})=\frac{\alpha}{n}\sum_{i,j}i^{|i|+|j|+|i||j|}\delta_{ij}=\alpha,
\end{equation}
as claimed. Meanwhile $C\ell_{p,q}\RR$ was discussed in section \ref{clifford}. Let $|x|_p=|x|-|x|_q\mod 2$. Then compute
\begin{align}\begin{split}
\cZ_{C\ell_{p,q}\RR}(\RR P^2_1)
&=\frac{1}{2^{(p+q)}}\sum_{N_i}\eta(\mathds{1}\otimes\tau)\left(\Gamma_1^{N_1}\cdots\Gamma_n^{N_n}\otimes *(\Gamma_1^{N_1}\cdots\Gamma_n^{N_n})\right)\\
&=\frac{1}{2^{(p+q)}}\sum_{N_i}i^{|x|}(-1)^{|x|_q}\eta\left(\Gamma_1^{N_1}\cdots\Gamma_n^{N_n}\otimes\Gamma_1^{N_1}\cdots\Gamma_n^{N_n}\right)\\
&=\frac{\alpha}{2^{(p+q)/2}}\sum_{N_i}i^{|x|}(-1)^{|x|_q}(-1)^{\{x\}(\{x\}-1)/2}\\
&=\frac{\alpha}{2^{(p+q)/2}}\sum_{N_i}i^{\{x\}}(-1)^{|x|_q}\\
&=\frac{\alpha}{2^{(p+q)/2}}\sum_{N_i}i^{|x|_p}(-i)^{|x|_q}\\
&=\alpha\left(\frac{1}{2^{p/2}}\sum_{k=0}^p\left(\begin{array}{c}p\\k\end{array}\right)i^k\right)\left(\frac{1}{2^{q/2}}\sum_{l=0}^q\left(\begin{array}{c}q\\l\end{array}\right)(-i)^{l}\right)\\
&=\alpha\exp((p-q)\pi i/4).
\end{split}\end{align}

This completes our argument. As a consistency check, one may evaluate the state sums on other closed pin manifolds and verify that they yield powers of the ABK invariant. This was done in ref. \cite{BT14} for orientable pin (spin) surfaces. They show that $C\ell_{1,0}\RR$ yields partition function $\cZ(M^\text{or})\sim\text{Arf}(M^\text{or})=\text{ABK}(M^\text{or})\in\{\pm 1\}$.

\subsubsection*{Acknowledgments}

The author is especially grateful to John Barrett for acquainting me with the key references \cite{BT14,P84}, for productive discussions in the early stages of this project, and for sharing his thoughts on the manuscript. I would also like to thank Dave Aasen, Matt Heydeman, Nick Hunter-Jones, Anton Kapustin, and Eric Samperton for enlightening conversations and Arun Debray for valuable commentary on the manuscript. This work was performed partly at Fudan University and was supported in part by the U.S. Department of Energy, Office of Science, Office of High Energy Physics, under Award Number de-sc0011632 as well as the Simons Investigator Award.

\pagebreak


\begin{thebibliography}{99}

\bibitem{K87} L. H. Kauffman, ``State models and the Jones polynomial,'' Topology, 26(3), 395-407 (1987).

\bibitem{SR16} K. Shiozaki, S. Ryu, ``Matrix Product States and Equivariant Topological Field Theories for Bosonic Symmetry-Protected Topological Phases in 1+1 Dimensions,'' arXiv:  1607.06504

\bibitem{SSGR17} K. Shiozaki, H. Shapourian, K. Gomi, S. Ryu, ``Many-Body Topological Invariants for Fermionic Short-Range Entangled Topological Phases Protected by Antiunitary Symmetries,'' arXiv: 1710.01886.

\bibitem{KTY16a} A. Kapustin, A. Turzillo, M. You, ``TQFT and MPS,'' Phys. Rev. B 96, 075125 (2017).

\bibitem{KTY16b} A. Kapustin, A. Turzillo, M. You, ``Spin TQFT and Fermionic MPS,'' Phys. Rev. B 98, 125101 (2018).

\bibitem{FK11} L. Fidkowski, A. Kitaev, ``Topological Phases of Fermions in One Dimension,'' Phys. Rev. B, 83, 075103 (2011).

\bibitem{BT14} J.~Barrett, S.~Tavares, ``Two-Dimensional State Sum Models and Spin Structures,'' Commun. Math. Phys. 336, 63-100 (2015).

\bibitem{NR15} S. Novak, I. Runkel, ``State Sum Construction of Two-Dimensional TQFTs on Spin Surfaces,'' J. Knot Theory Ramifications 24 no.05, 1550028 (2015).

\bibitem{DG18} A. Debray, S. Gunningham, ``The Arf-Brown TQFT of $\text{Pin}^-$ Surfaces,'' arXiv: 1803.11183.

\bibitem{TY17} A. Turzillo, M. You, ``Fermionic MPS and 1D SRE Phases with Anti-Unitary Symmetries,'' Phys. Rev. B 99, 035103 (2019).

\bibitem{BWHV17} N. Bultinck, D. J. Williamson, J. Haegeman, F. Verstraete, ``Fermionic matrix product states and one-dimensional topological phases,'' Phys. Rev. B95,075108 (2017).

\bibitem{LP06} A. Lauda, H. Pfeiffer, ``State sum construction of two-dimensional open-closed Topological Quantum Field Theories,'' arXiv: math/0602047

\bibitem{ABS64} M. Atiyah, R. Bott, A. Shapiro, ``Clifford Modules,'' Topology, 3 (Suppl. 1): 3–38 (1964).

\bibitem{KT89} R. C.~Kirby, L. R.~Taylor, ``Pin structures on low dimensional manifolds,'' Geometry of low-dimensional manifolds, 2, 177-242 (1989).

\bibitem{P84} U.~Pinkall, ``Regular Homotopy Classes of Immersed Surfaces,'' Topology Vol. 24. No. 4, 421-434 (1985).

\bibitem{JT66} I. James, E. Tomas, ``Note on the classification of cross-sections,'' Topology. 4 351-359 (1966).

\bibitem{S77} K. Schlichting, ``Regultire Homotopie geschlossener Flachen,'' Diplomarbeit TU Berlin (1977). 

\bibitem{MP78} W. H. Meeks III, J. Patrusky, ``Representing homology classes by embedded circles on a compact surface,'' Ill. J. Math. 22 262-269 (1978).

\bibitem{B72} E.H.~Brown, ``Generalizations of the Kervaire invariant,'' Ann. Math. 95 368-383 (1972).

\bibitem{S04} G.~Segal, ``The Definition of a Conformal Field Theory,” in Topology, Geometry and Quantum Field Theory, ed. U. Tillmann, London Math. Soc. Lect. Note Ser., Vol. 30 8, p.421-577 (2004).

\bibitem{MS06} G. W. Moore, G. Segal, ``D-branes and K-theory in 2D topological field theory,'' arXiv: hep-th/0609042.

\bibitem{TT06} V.~Turaev, P.~Turner, ``Unoriented topological quantum field theory and link homology,'' Algebraic \& Geometric Topology 6, 1069 (2006).

\bibitem{FHK94} S. Fukuma, A. Hosono, H. Kawai, ``Lattice Topological Field Theory in two dimensions,'' Comm. Math. Phys. 161, no. 1, 157-175 (1994).

\bibitem{K90} L. H. Kauffman, ``An invariant of regular isotopy,'' Trans. Am. Math. Soc. 318, 417-471 (1990).

\bibitem{FY89} P. J. Freyd, D. N. Yetter, ``Braided compact closed categories with applications to low-dimensional topology,'' Adv. Math.,  77(2) 156-182 (1989).

\bibitem{Y89} D. N. Yetter. ``Category theoretic representations of knotted graphs in s3,'' Adv. Math., 77(2) 137-155 (1989).

\bibitem{K05} L. H. Kauffman. ``Knot diagrammatics,'' Handbook of Knot Theory, 233-318, Elsevier (2005).

\bibitem{RT90} N. Yu. Reshetikhin, V. G. Turaev, ``Ribbon graphs and their invariants derived from quantum groups,'' Comm. Math. Phys., 127(1) 1–26 (1990).

\bibitem{V07} O. Viro, ``Quantum Relatives of the Alexander Polynomial,'' St. Petersburg Math. J., Vol. 18, No. 3, 391-457 (2007).

\bibitem{P91} U. Pachner, ``P.l. homeomorphic manifolds are equivalent by elementary shellings,'' Eur. J. Comb., 12(2) 129-145 (1991).

\bibitem{L99} W. B. R. Lickorish, ``Simplicial moves on complexes and manifolds,'' In Proceedings of the Kirbyfest (Berkeley, CA, 1998), volume 2 of Geom. Topol. Monogr., 299-320 Geom. Topol. Publ., Coventry (1999).

\bibitem{CL14} A. Coward, M. Lackenby, ``An upper bound on Reidemeister moves,'' American Journal of Mathematics, 136 1023-1066 (2014).

\bibitem{T10} V. Turaev, ``Homotopy Quantum Field Theory,'' EMS Tracts in Mathematics Vol. 10 (2010).

\bibitem{FH16} D. S. Freed, M. J. Hopkins, ``Reflection positivity and invertible topological phases,'' arXiv: 1604.06527.

\bibitem{D99} P. Deligne, ``Notes on Spinors,'' Quantum Fields and Strings: a Course for Mathematicians, vol. 1, AMS. Providence, RI, pp. 99–135 (1999).

\bibitem{W64} C. T. C. Wall, ``Graded Brauer Groups,'' Journal f\"ur die reine und angewandte Mathematik, 213, 187-199 (1964).

\bibitem{L05} T.-Y. Lam, ``Introduction to Quadratic Forms over Fields,'' Graduate Studies in Mathematics Vol. 67, American Mathematical Soc. (2005).

\bibitem{A41} C. Arf, ``Untersuchungen \"uber quadratische Formen in K\"orpern der Charakteristik 2, I,'' J. Reine Angew. Math. 183, 148-167 (1941).

\bibitem{J80} D. Johnson, ``Spin structures and quadratic forms on surfaces,'' J. London Math. Soc. (2) 22 no. 2, 365-373 (1980).

\bibitem{Y18} K. Yonekura, ``On the cobordism classification of symmetry protected topological phases,'' arXiv: 1803.10796.

\bibitem{FM04} D. S. Freed, G. W. Moore, ``Setting the quantum integrand of M-theory,'' Commun. Math. Phys. 263, 89-132 (2006).

\end{thebibliography}
\end{document}